\documentclass[11pt]{article}
\usepackage{geometry}                
\usepackage{graphicx}
\usepackage{amssymb}
\usepackage{epstopdf}
\begin{document}

\title{{\bf Solutions of the cubic Fermat equation in ring class fields of imaginary quadratic fields (as periodic points of a 3-adic algebraic function)} }        
\author{Patrick Morton}        
\date{}          
\maketitle

\begin{abstract}
Explicit solutions of the cubic Fermat equation are constructed in ring class fields $\Omega_f$, with conductor $f$ prime to $3$, of any imaginary quadratic field $K$ whose discriminant satisfies $d_K \equiv 1$ (mod $3$), in terms of the Dedekind $\eta$-function.  As $K$ and $f$ vary, the set of coordinates of all solutions is shown to be the exact set of periodic points of a single algebraic function and its inverse defined on natural subsets of the maximal unramified, algebraic extension $\textsf{K}_3$ of the $3$-adic field $\mathbb{Q}_3$.  This is used to give a dynamical proof of a class number relation of Deuring.  These solutions are then used to give an unconditional proof of part of Aigner's conjecture: the cubic Fermat equation has a nontrivial solution in $K=\mathbb{Q}(\sqrt{-d})$ if $d_K \equiv 1$ (mod $3$) and the class number $h(K)$ is not divisible by $3$.  If $3 \mid h(K)$, congruence conditions for the trace of specific elements of $\Omega_f$ are exhibited which imply the existence of a point of infinite order in $Fer_3(K)$.
\end{abstract}

\section{Introduction.}  

In the 1950's Aigner wrote a series of papers [1]-[4] on the cubic Fermat equation in quadratic fields, building on work of Fueter [18].  Most of Aigner's results had to do with the nonexistence of nontrivial solutions to the cubic Fermat equation in quadratic fields.  The character of solutions to this equation in $K=\mathbb{Q}(\sqrt{m})$ is the same as in the companion field $K_1=\mathbb{Q}(\sqrt{-3m})$: if there is a nontrivial solution in one field, then the same holds for the companion field.  Thus, there are four families of coupled quadratic fields to consider: the fields generated over $\mathbb{Q}$ by

$$1: \sqrt{3n+1} \quad \textrm{or} \quad \sqrt{-(9n+3)},$$
$$2: \sqrt{-(3n+1)} \quad \textrm{or} \quad \sqrt{9n+3},$$
$$3: \sqrt{3n+2} \quad \textrm{or} \quad \sqrt{-(9n+6)},$$
$$4: \sqrt{-(3n+2)} \quad \textrm{or} \quad \sqrt{9n+6},$$ \smallskip

\noindent where $n \ge 0$.  Fueter had proved in 1913 that there is no nontrivial solution of $x^3+y^3=z^3$ in quadratic fields of the second family, when $m=-(3n+1)$ and $3$ does not divide the class number $h(K)$ of the imaginary field $K=\mathbb{Q}(\sqrt{m})$ (see [18]).   Aigner [1] proved impossibility of solutions for the first family when the class number of the imaginary field is not divisible by $3$.  Aigner's further investigations (see [35]) led him to conjecture in [3, p.16] that there always exist nontrivial solutions in the fields of the third and fourth families.  The recent paper of Jones and Rouse [25] contains strong evidence for the truth of Aigner's conjecture.  They prove a conditional result, based on the Birch-Swinnerton-Dyer conjecture, for the existence of nontrivial solutions to the cubic Fermat equation in a quadratic field, in terms of the number of solutions of certain ternary quadratic forms; and their result, if shown to be unconditionally true, would imply Aigner's conjecture.  \medskip

In this paper we study solutions of the equation
$$Fer_3: \hspace{.1 in} 27X^3+27Y^3=X^3Y^3$$
in ring class fields of the imaginary quadratic fields of the {\it fourth} family:
$$K=\mathbb{Q}(\sqrt{-d}), \ \ -d \equiv 1 \ (\textrm{mod} \ 3),$$
and show that they arise as periodic points of a specific (multi-valued) algebraic function $T(z)$ over $\mathbb{C}$.  Moreover, it turns out that these solutions give all the complex periodic points of $T(z)$ other than $z=3$, and can be given explicitly in terms of modular functions.  \medskip

Recall that the ring class field $\Omega_f$ of an imaginary quadratic field $K=\mathbb{Q}(\sqrt{-d})$ is the unique finite algebraic extension of $K$ whose Galois group $\textrm{Gal}(\Omega_f/K)$ is abelian, having the property that the prime ideals (relatively prime to $f$) of the ring of integers $R_K$ of $K$, which split completely into prime ideals of degree $1$ in the ring of integers $R_{\Omega_f}$ of $\Omega_f$, are exactly the prime ideals $\mathfrak{p}$ for which $\mathfrak{p}=(\xi)$ is principal in $R_K$ with $\xi \equiv r$ (mod $f$) and $r \in \mathbb{Z}$.  It then follows that $\textrm{Gal}(\Omega_f/K) \cong A_f/P_f$, where $A_f$ is the group of fractional ideals of $K$ which are relatively prime to $f$ and $P_f$ is the subgroup of $A_f$ consisting of principal ideals of the form $(\xi)$ for numbers $\xi \equiv r$ (mod $f$) and $r \in \mathbb{Z}$.   The set of all such integers $\xi$ of $R_K$ is a ring, which gives rise to the name {\it ring class field}.  The properties of these fields are developed in the classical theory of complex multiplication.  (See [10], [15], and the paper of Hasse [21].)  \medskip

In Sections 2 and 3 we shall exhibit an explicit solution of $Fer_3$ in the ring class field $\Omega_f$ of $K$ whose conductor $f$ is relatively prime to $3$, in terms of the Dedekind eta function $\eta(z)$, defined as 
$$\eta(z)=q^{1/24} {\prod_{n=1}^\infty}\left(1-q^n\right),\hspace{.1 in} q = e^{2\pi i z}, \ \ z \in \mathbb{H},$$ where $\mathbb{H}$ is the upper half-plane.  \bigskip

\noindent {\bf Theorem 1.} {\it If $K=\mathbb{Q}(\sqrt{-d})$, with $-d =d_Kf^2 \equiv 1$ (mod 3) and $d_K=\textrm{disc}(K/\mathbb{Q})$, let
$$w=\cases{k+\frac{\sqrt{-d}}{2}, &if $2 \mid d$,\cr \frac{k+\sqrt{-d}}{2}, &if $(2,d)=1$;} \eqno{(1.1)}$$
where $k^2 \equiv -d/4$ resp. $-d$ (mod 9) and $k \equiv 1$ (mod 6).  Then a non-trivial solution $(X,Y)=(\alpha,\beta)$ of $Fer_3$ in the ring class field $\Omega_f$ of $K$ is given by
$$\alpha =3+\gamma^3, \ \ \beta=\frac{3(\alpha^\tau+6)}{\alpha^\tau-3},$$ \smallskip
where $\gamma$ is the unique multiple
$$\gamma=\omega^i \frac{\eta(w/9)}{\eta(w)}, \ \ \ \omega=\frac{-1+\sqrt{-3}}{2}, \ i \in \{0,1,2\},$$ \smallskip
of $\eta(w/9)/\eta(w)$ which lies in $\Omega_f$, and $\tau=(\Omega_f/K, \wp_3) \in \textrm{Gal}(\Omega_f/K)$ is the Artin symbol for the prime ideal $\wp_3 = (3,w)$ of $R_K$.  The numbers $\alpha$ and $\beta$ are conjugate algebraic integers over $\mathbb{Q}$, and $\Omega_f=\mathbb{Q}(\alpha)=\mathbb{Q}(\gamma)$.  In addition, $(\gamma)=\gamma R_{\Omega_f}=\wp_3'R_{\Omega_f}$, where $\wp_3'$ is the conjugate ideal of $\wp_3$ in $R_K$. } \bigskip

Parts of this theorem follow from known results.  The fact that $\alpha, \beta \in \Omega_f$ follows easily from a theorem of Schertz [36, p. 159].  The difficult part of the theorem is showing that $\gamma \in \Omega_f$ for some $i$.  In the case $f=1$, this follows from a result of Kubert and Lang [26, p. 296], who extended an old result of Fricke, Hasse, and Deuring.  However, the case $f>1$ seems not to have been noted before (see the Remark following Theorem 3.4).  The numbers $\alpha$ and $\beta$ turn out to be quotients of values of cubic theta functions (see \S 2, as well as [31] and the references in that paper). \medskip

For the sake of completeness we give a proof of Theorem 1 for general $f$ (relatively prime to $3$), using only the theory of complex multiplication and a modular parametrization of $Fer_3$ that was worked out in [31].  Our discussion shows that the point $(\alpha,\beta) \in Fer_3(\Omega_f)$ has infinite order, and corresponds to an integral point $P_d=(\frac{\alpha^\tau+6}{\alpha},\alpha^\tau+6)$ on the isomorphic elliptic curve $E: Y^2-9Y=X^3-27$.  (See Theorem 6.6.) \medskip

If $h(-d)=h(\textsf{R}_{-d})=[\Omega_f:K]$ denotes the class number of the order $\textsf{R}_{-d}$ of discriminant $-d$ in $K$, then the minimal polynomial $p_d(x)$ of $\alpha$ over $\mathbb{Q}$ has degree $2h(-d)$ and satisfies the transformation formula
$$(x-3)^{2h(-d)} p_d\bigg(\frac{3(x+6)}{x-3}\bigg)=3^{3h(-d)}p_d(x),$$  \smallskip
while the minimal polynomial $q_d(x)$ of $\gamma$ over $\mathbb{Q}$ satisfies
$$x^{2h(-d)}q_d\bigg(\frac{3}{x}\bigg)=3^{h(-d)} q_d(x).$$ \smallskip
Furthermore, the roots of $p_d(x)$ are periodic points of the dynamical system defined by a single algebraic function and its inverse function in the maximal unramified, algebraic extension $\textsf{K}_3$ of the $3$-adic field $\mathbb{Q}_3$.  \bigskip

\noindent {\bf Theorem 2.} {\it The roots of the polynomials $p_d(x)$ (as $-d$ varies over quadratic discriminants $\equiv 1$ (mod $3$)) which are prime to the respective ideals $\wp_3=(3,w)$ in $R_K$ are {\it all} the periodic points of the algebraic function
$$T(z)=\frac{z^2}{3}(z^3-27)^{1/3}+\frac{z}{3}(z^3-27)^{2/3}+\frac{z^3}{3}-6\eqno{(1.2)}$$
in its domain $\{z: |z|_3 \ge 1\} \subset \textsf{K}_3$, under the natural embedding $L \rightarrow L_\mathfrak{p} \subset \textsf{K}_3$, where $L=\Omega_f$ and $L_\mathfrak{p}$ is the completion of $L$ at a prime ideal $\mathfrak{p}$ for which $\mathfrak{p} \vert \wp_3$.  The remaining roots of the polynomials $p_d(x)$ (for the same $d$'s) are, together with $z=3$, all the periodic points of the inverse algebraic function
$$S(z)=\frac{z+6}{(z^2+3z+9)^{1/3}}$$ \smallskip
in the disk $\textsf{D}_3=\{z: |z-3|_3 \le |3|_3^3 \} \subset \textsf{K}_3$.} \bigskip

Viewing $T(z)$ and $S(z)$ as single-valued functions defined on the subsets $|z|_3 \ge 1$ and $\textsf{D}_3$ of the unramified extension $\textsf{K}_3$, respectively, is one way of choosing specific branches of the algebraic functions $T$ and $S$, since the cube roots of unity are not contained in $\textsf{K}_3$; and this is what makes it possible to consider their iterates on these sets.  Note also that $T$ and $S$ are conjugate maps on their respective domains (see equation (4.20)).  \medskip

The analysis of the maps $T(z)$ and $S(z)$ has the following consequence.  Let $\mathfrak{D}_n$ denote the set of discriminants $-d \equiv 1$ (mod 3) of orders in imaginary quadratic fields $K=\mathbb{Q}(\sqrt{-d})$ for which the Frobenius automorphism $\tau=(\Omega_f/K, \wp_3)$ in the corresponding ring class field $\Omega_f$ has order $n$.  Then 
$$\sum_{-d \in \mathfrak{D}_n}{h(-d)}=nN_3(n)=\sum_{k|n}{\mu (n/k)3^k}, \quad n>1.\eqno{(1.3)}$$ \medskip
This equivalent to a special case of a class number formula of Deuring [12, p. 269] (corrected in [13, p. 24]; see Section 4), which he derived as a consequence of his theory of the endomorphism rings of elliptic curves.  Here we give an alternative, {\it dynamical} interpretation of (1.3) by deriving this formula as a consequence of the arithmetic of the cubic Fermat equation and the dynamical system defined by (1.2) on $\textsf{K}_3$.  Then the number on either side of (1.3) is the number of periodic points $\alpha$ of $T(z)$ of minimal period $n$, all of which lie in the unit group $\textsf{U}_3$ of the ring of $3$-adic integers of $\textsf{K}_3$.  \medskip

Thus, the function $T(z)$ locates the {\it global} numbers $\alpha$ and $\beta$ which are the roots of $p_d(x)$ inside of the unramified {\it local} $3$-adic field $\textsf{K}_3$.  This is analogous to the way in which we normally envision algebraic numbers over $\mathbb{Q}$ as elements of the complex field $\mathbb{C}$.  For the infinite set of algebraic numbers defined by the polynomials $p_d(x)$, the map $T$ is a common lifting of the Frobenius map to the $3$-adic field $\textsf{K}_3$, which is generated over $\mathbb{Q}_3$ by these numbers.  Note the connection to the fact that the unique unramified extension of degree $n$ over $\mathbb{Q}_3$  is generated by primitive roots of unity of order $3^n-1$, which are periodic points of the map $F(z) = z^3$ with minimal period $n$.  Equation (1.3) is an expression of the fact that $T(z)$ and $F(z)$ have the same number of periodic points of minimal period $n$ in the unit group $\textsf{U}_3$, for all $n>1$.  \medskip

The $3$-adic proof of Theorem 2 actually yields more, that the only periodic points, suitably defined, of the multi-valued function $T(z)$ on either the algebraic closure $\overline{\mathbb{Q}}_3$ of $\mathbb{Q}_3$ or the complex field $\mathbb{C}$ are $z=3$ and the roots of the polynomials $p_d(x)$.  (See Theorem 4.4.)  Thus, the periodic points $z \neq 3$ of $T(z)$ in $\mathbb{C}$ have number theoretic significance, in that they generate ring class fields, unramified over the prime $p=3$, of the imaginary quadratic fields $K$ considered in Theorems 1 and 2.  It is also possible to show that all {\it pre-periodic} points of $T(z)$ other than $z=3\omega, 3\omega^2$ generate ring class fields, {\it ramified} over $p=3$, of the same quadratic fields considered above. (See Section 5.)  \medskip

In another paper [32] we will show that similar considerations apply to the solutions of the quartic Fermat equation which are studied in [28].  Their coordinates also represent periodic points of an algebraic function, which is defined on a subset of the maximal unramified, algebraic extension $\textsf{K}_2$ of the $2$-adic field $\mathbb{Q}_2$, and which yields the analogue of the dynamical relation (1.3) for the prime $p=2$.  Based on these examples, and the fact that Deuring's class number formulas hold for all primes $p$, it is reasonable to make the following conjecture. (This conjecture is also stated in [32].)  \medskip

To state the conjecture, define an imaginary quadratic field $K$ to be $p$-{\it admissible}, for a given prime $p \in \mathbb{Z}$, if $\displaystyle \left(\frac{d_K}{p}\right)=+1$, where $d_K$ is the discriminant of $K$, so that $p$ splits into two prime ideals in the ring of integers $R_K$ of $K$. \bigskip

\noindent {\bf Conjecture.} {\it Let $p$ be a fixed prime number.  There is an algebraic function $T_p(z)$, defined on a certain subset of the maximal unramified, algebraic extension $\textsf{K}_p$ of the $p$-adic field $\mathbb{Q}_p$, with the following properties: \smallskip

1)  All ring class fields of any $p$-admissible quadratic field $K \subset \mathbb{Q}_p$ whose conductors are relatively prime to $p$ are generated over $K$ by periodic points of $T_p(z)$ contained in the unramified extension $\textsf{K}_p$.  \smallskip

2) All ring class fields of $K$ whose conductors are divisible by $p$ are generated over $K$ by pre-periodic points of $T_p(z)$ contained in the algebraic closure $\overline{\mathbb{Q}}_p$. \smallskip

3) All but finitely many of the periodic and pre-periodic points of $T_p(z)$ contained in $\overline{\mathbb{Q}}_p$ generate ring class fields over suitable $p$-admissible quadratic fields $K$.}  \bigskip

The results of this paper show that this conjecture is true for the prime $p=3$, and further papers will show that it is also true for $p=2$ and $p=5$. \medskip

We finish this paper by proving Aigner's conjecture for any imaginary quadratic field $K=\mathbb{Q}(\sqrt{-d})$ of the fourth family whose class number $h(K)$ is not divisible by 3, using the solutions $(\alpha,\beta)$ from Theorem 1 and elementary facts about elliptic curves.  In particular, Aigner's conjecture holds unconditionally for a set of quadratic discriminants $d_K \equiv 1$ (mod $3$) which has density at least $1/2$ in the set of all such discriminants.  Similar arguments (see Theorem 6.4) show that for any discriminant $-d \equiv 1$ (mod $3$) there are infinitely many positive fundamental discriminants $D$ for which the quartic field $L=\mathbb{Q}(\sqrt{-d},\sqrt{D})$ contains a nontrivial solution of $Fer_3$.  This holds for any positive, square-free $D \equiv 1$ (mod $12$) for which $(d,D)=1$ and $h(-dD) \not \equiv 0$ (mod $3$).  \medskip

We also prove that Aigner's conjecture is true for $K$ when $3 \mid h(K)$, if the solution $(\alpha,\beta)$ of $Fer_3$ constructed in Theorem 1 for the Hilbert class field $\Sigma=\Omega_1$ satisfies
$$Tr_{\Sigma/K}(1/\alpha) \not \equiv 0 \ (\textrm{mod} \ \wp_3^2).\eqno{(1.4)}$$
This condition, which follows from a formal group calculation in Section 7, is equivalent to the condition that $9 \nmid H_{-d}'(6)$, where $H_{-d}'(X)$ is the derivative of the class equation for the discriminant $-d$.  With this criterion, we have verified that Aigner's conjecture holds for all $39$ fields of the fourth family for which the class number satisfies $h(K)=3,6,9$, or $12$; and (1.4) holds for $38$ of these $39$ fields (see Propositions 7.3 and 7.4).  In the one remaining case, $d=2132$, a similar condition modulo $\wp_3^3$ shows the truth of Aigner's conjecture for the field $K=\mathbb{Q}(\sqrt{-2132})$ (Theorem 7.5). \medskip

These results fit very nicely with the results of Fueter and Aigner.  In their case, the condition $3 \nmid h(K)$ for the class number of an imaginary quadratic field $K$ in families 1 or 2 is used to prove {\it nonexistence} of solutions to $Fer_3$ in $K$ and its companion field, while here the corresponding condition for the fourth family guarantees the {\it existence} of solutions in $K$.  Aigner showed, in addition, that various congruence conditions, such as $2$ being a cubic non-residue of the prime factors $p \equiv 1$ (mod $3$) of $m$, or $2^{1/3}-1$ being a quadratic non-residue of a prime factor $p \equiv 5$ (mod $6$) of $m$, is sufficient in some cases to prove nonexistence for fields $K=\mathbb{Q}(\sqrt{m})$ in those same families.  See [1, p. 249], [2, p. 335], and [35, pp. 283-285].  The congruence condition (1.4), on the other hand, guarantees the {\it existence} of a solution in $K$.  In our case, any nontrivial point $Q \in Fer_3(K)$, i.e. one which is not in $Fer_3(\mathbb{Q})$, has infinite order on $Fer_3$ [35, p. 281], so (1.4) guarantees that the rank of $Fer_3(K)$ is positive.  See Theorem 7.1 and its Corollary 1 for a simple proof using the formal group of the curve $E$ mentioned above.

\section{Constructing solutions of $Fer_3$ in ring class fields.}

The impetus for discovering the solutions in Theorem 1 came from studying a relationship between between $Fer_3$ and the elliptic curve
$$E_\alpha: \hspace{.1 in} Y^2+\alpha XY + Y = X^3 \eqno{(2.1)}$$ \smallskip
\noindent in Deuring normal form (with base point $O$).  The curve $Fer_3$ parametrizes the 3-isogenies with kernel
$$T_3=\{O,(0,0),(0,-1)\}\subset E_\alpha$$
from one elliptic curve in Deuring normal form to another.  If $(\alpha,\beta) \in Fer_3(k)$, where the characteristic of  $k$ is not 3, then there is an isogeny $\phi_{\alpha,\beta}: E_\alpha \rightarrow E_\beta$  with kernel $T_3$ which is defined over $k(\alpha,\beta,\sqrt{-3})$.  Further, let $G_{12}$ be the group of linear fractional mappings in $z$ generated by
$$\sigma_1(z)=\frac{3(z+6)}{z-3},\hspace{.1 in} \sigma_2(z)=\omega z, \hspace{.1 in} \omega^2+\omega+1=0.$$ \smallskip
\noindent If $\phi:E_\alpha \rightarrow E_\beta$ is any isogeny with kernel $T_3$, then for some $\sigma \in G_{12}$, the point $(\sigma(\alpha),\beta)$ lies on $Fer_3$.  See [29] for proofs of these facts.  In addition, the remaining points of order 3 on $E_\alpha$ can be given in terms of simple rational functions of points $(\alpha,\beta)$ on $Fer_3$: if $\alpha \neq 0$, the points of order 3 on $E_\alpha$ are $(0,0)$, $(0,-1)$, and the points
$$(x,y)=\left(\frac{-3\beta}{\alpha(\beta-3)},\frac{\beta - 3 \omega^i}{\beta-3}\right),\hspace{.1 in} i = 1,2,\eqno{(2.2)}$$ \smallskip
\noindent where $\beta$ runs over the three elements of $\bar k$ for which $(\alpha,\beta)$ is a point on $Fer_3$.  Exchanging $\alpha$ and $\beta$ in this formula gives the points of order 3 on the isogenous curve $E_\beta$.  \bigskip

The points $(X,Y)$ on $Fer_3$ can be parametrized by the following functions on the upper half-plane $\mathbb{H}$:
$$X=\mathfrak{f}(z)=3+\left(\frac{\eta(z/9)}{\eta(z)}\right)^3, \ \ Y= \mathfrak{g}(z)=3+\left(\frac{3\eta(3z)}{\eta(z/3)}\right)^3, \ \ z \in \mathbb{H},\eqno{(2.3)}$$
which are both modular functions for the group $\Gamma(9)$.  (See [31].  Note that $\tau$ has been replaced by $z/3$ in the formulas of [31, Thm. 12].)  In fact, $\mathfrak{f}(z)$ and $\mathfrak{g}(z)$ generate the field $\textsf{K}_{\Gamma'}$ of modular functions for the subgroup
$$\Gamma' = \{\left(\begin{array}{cc}a & 9b \\3c & d \end{array}\right) \in \Gamma: a\equiv d \equiv \pm 1 \ (\textrm{mod} \ 3), b, c \in \mathbb{Z} \}$$
of $\Gamma=\textrm{SL}(2,\mathbb{Z})$. (This follows from the remark on p. 362 of [31].)  If $R=\left(\begin{array}{cc} 9 & 0 \\9 & 1 \end{array}\right)$, then $R^{-1} \Gamma' R = \Gamma_0(27)$.  This implies that $\textsf{K}_{\Gamma'}$ is isomorphic to the function field of the curve $X_0(27)$, by the map $f(z) \rightarrow f(Rz)$, and therefore $X_0(27) \cong Fer_3$, a fact which is known.  The parametrization (2.3) is related to cubic theta functions, which were introduced in 1991 by the Borweins.  (See [31, \S 5] and the references in that paper.)  The function $\mathfrak{f}(z)-3$ occurs as the function $y(w)$ in Weber [39, p.255].  \medskip

We will take $z$ equal to the number $w$ defined in Theorem 1, so that we have:
$$\alpha=\mathfrak{f}(w), \ \ \ \beta =\omega^i \mathfrak{g}(w), \ \ \omega=\frac{-1+\sqrt{-3}}{2}, \ \ i \in \{0,1,2\}.\eqno{(2.4)}$$
(We will see below that $i \neq 0$.)  For $z=w$ it follows immediately from a theorem of Schertz [36, Thm. 6.6.4, p. 159] that $\alpha=\mathfrak{f}(w)$ lies in the ring class field $\Omega_f$, where $-d=d_Kf^2,$ using the notation of Theorem 1.  This is because $\{9, w\}$ is a basis for the ideal $\wp_{3,-d}^2=(3,w)^2$ of the order $\textsf{R}_{-d} \subset K=\mathbb{Q}(\sqrt{-d})$.  Since $\{1, 3w\}$ is a basis for a proper ideal of $\textsf{R}_{-9d}$, the order of conductor $3f$, Schertz's theorem also implies that $\displaystyle \left(\frac{\eta(3w)}{\eta(w/3)}\right)^3 \in \Omega_{3f}=\Omega_f(\omega)$.  However, $(\alpha,\beta)$ is a point on $Fer_3$, and this implies that $\displaystyle \beta^3=\frac{27\alpha^3}{\alpha^3-27} \in \Omega_f$.  Since $[\Omega_{3f}: \Omega_f]=2$, $x^3-\beta^3$ must have a linear factor over $\Omega_f$, so one of the cube roots of $\beta^3$ must lie in $\Omega_f$.  By choosing $i$ appropriately in (2.4), we therefore have that $\beta \in \Omega_f$.  This proves   \bigskip

\noindent {\bf Proposition 2.1.} {\it Let $(\alpha, \beta)$ be the point on $Fer_3$ given by (2.4), where $w \in \textsf{R}_{-d}\subset K=\mathbb{Q}(\sqrt{-d})$ is defined by (1.1), and $-d=d_Kf^2$.  For some $i \in \{0,1,2\}$ the numbers $\alpha, \beta$ lie in the ring class field $\Omega_f$ of conductor $f$ over $K$.}  \bigskip

Note that the $j$-invariant of the elliptic curve $E_\alpha$ in (2.1) is
$$j(E_\alpha) = j_\alpha=\frac{\alpha^3(\alpha^3-24)^3}{\alpha^3-27};\eqno{(2.5)}$$
and that $\alpha \neq 3\omega^i$ ($0 \le i \le 2$).  This is easy to see, since
$$\left(\frac{3}{\alpha}\right)^3+\left(\frac{3}{\beta}\right)^3=1,$$
and $\beta$ is a value of the holomorphic function $\mathfrak{g}(z)$ on $\mathbb{H}$.  In fact, neither of $\mathfrak{f}(z)$ and $\mathfrak{g}(z)$ is infinite on $\mathbb{H}$, so neither can take the value $3\omega^i$ for $z \in \mathbb{H}$. \medskip

From [31, \S 4] we take the fact that for
$$(u,v) \in \{(0,1/3), (1/3, 0), (1/3, 1/3), (1/3, -1/3)\},$$
the four functions $A_{(u,v)}(z)=l_{(u,v)}(z)^3$, with
$$l_{(u,v)}(z)=\cases{3+\left(\frac{\eta \left(\frac{uz+v}{3}\right)}{\eta(z)}\right)^3, \ &for  $(u,v) \neq (0,1/3)$, \cr 3+27\left(\frac{\eta(9z)}{\eta(z)}\right)^3, \ &for $(u,v)=(0,1/3),$}$$
are the roots of the equation
$$X(X-24)^3-j(z)(X-27)=0.\eqno{(2.6)}$$
Note that $l_{(1/3,0)}(z)=\mathfrak{f}(z)$ and $l_{(0,1/3)}(z/3)=\mathfrak{g}(z)$ from (2.3). Thus, (2.6) gives that
$$j(z)=\frac{\mathfrak{f}(z)^3(\mathfrak{f}(z)^3-24)^3}{\mathfrak{f}(z)^3-27}=\frac{\mathfrak{g}(z)^3(\mathfrak{g}(z)^3+216)^3}{(\mathfrak{g}(z)^3-27)^3},\eqno{(2.7a)}$$
where the second equality follows on replacing $\mathfrak{f}(z)^3$ by $27\mathfrak{g}(z)^3/(\mathfrak{g}(z)^3-27)$ in the first formula.  We will also need the companion formula
$$j(z/3)=\frac{\mathfrak{f}(z)^3(\mathfrak{f}(z)^3+216)^3}{(\mathfrak{f}(z)^3-27)^3}=\frac{\mathfrak{g}(z)^3(\mathfrak{g}(z)^3-24)^3}{\mathfrak{g}(z)^3-27}.\eqno{(2.7b)}$$
Note that (2.7b) also follows from (2.7a) on putting $z/3$ for $z$ and using the fact that
$$\sigma_1(\mathfrak{f}(z))=\frac{3(\mathfrak{f}(z)+6)}{\mathfrak{f}(z)-3}=3+\left(\frac{3\eta(z)}{\eta(z/9)}\right)^3=\mathfrak{g}(z/3).\eqno{(2.8)}$$
This identity implies that
$$\frac{\mathfrak{g}(z/3)^3(\mathfrak{g}(z/3)^3+216)^3}{(\mathfrak{g}(z/3)^3-27)^3}=\frac{\mathfrak{f}(z)^3(\mathfrak{f}(z)^3+216)^3}{(\mathfrak{f}(z)^3-27)^3},$$
since the rational function $r(z)=z^3(z^3+216)^3/(z^3-27)^3$ is invariant under the linear fractional maps $(z \rightarrow \sigma(z))$ in the group $G_{12}$. (See the proof of Proposition 3.1.) \medskip

Now (2.4) and (2.7ab) imply that $x=\alpha=\mathfrak{f}(w)$ is a solution of 
$$x^3(x^3-24)^3-j(w)(x^3-27)=0,\eqno{(2.9)}$$
and that
$$j(w)=\frac{\alpha^3(\alpha^3-24)^3}{\alpha^3-27}, \hspace{.1 in} j\left(\frac{w}{3}\right)=\frac{\alpha^3(\alpha^3+216)^3}{(\alpha^3-27)^3}.\eqno{(2.10)}$$ \smallskip
Since $\{1,w\}$ is a basis for $\textsf{R}_{-d}$, we know that $j(w)$ is a root of $H_{-d}(x)$.  This shows that the curve $E_\alpha$, whose $j$-invariant is $j(w)$, has complex multiplication by $\textsf{R}_{-d}$.  Since $-d \equiv 1$ (mod 3), $\{3,w\}$ is a basis for a prime ideal divisor $\wp_{3,-d}=\wp_3 \cap \textsf{R}_{-d}$ of 3 in the order $\textsf{R}_{-d}$.  Thus, $j(w/3)$ is also a root of $H_{-d}(x)$ and $j(w)$ and $j(w/3)$ are conjugate over $\mathbb{Q}$.  By the theory of complex multiplication, both are generators of $\Omega_f$ over $K$.  \medskip

For later use we note that (2.4) and (2.9) also hold when $w=w_1/w_2$ and $\{w_1, w_2\}$ is a basis for an ideal $\mathfrak{a}$ in $\textsf{R}_{-d}$ such that $\{w_1,3w_2\}$ is a basis for $\wp_{3,-d} \mathfrak{a}$. \medskip

This discussion implies the following assertion.  \bigskip

\noindent {\bf Proposition 2.2.} {\it If $\alpha$ and $\beta$ in $\Omega_f$ are given by (2.4), the curves $E_\alpha$ and $E_\beta$ are isogenous, and the isogeny $\phi_{\alpha,\beta}=(E_\alpha \rightarrow E_\beta)$ represents a Heegner point on the curve $X_0(3)$.} \medskip

\noindent {\it Proof.} Formulas (2.4) and (2.7b) give that
$$j\left(\frac{w}{3}\right)=\frac{\beta^3(\beta^3-24)^3}{\beta^3-27}=j(E_\beta).$$
Hence, $j(E_\beta)$ is also a root of $H_{-d}(X)$, so that both $E_\alpha$ and $E_\beta$ have complex multiplication by $\textsf{R}_{-d}$.  This shows that the $3$-isogeny $\phi_{\alpha,\beta}$ satisfies the definition of a Heegner point.  (See [5], [6], [7], [19].)  $\square$ \bigskip

In terms of the notation in the paper of Gross [19], the isogeny $(E_\alpha \rightarrow E_\beta)$ in this proposition is the dual of the isogeny $(\mathcal{O},\mathfrak{n},[\mathfrak{a}])=(\textsf{R}_{-d},\wp_3,[\wp_3])$.  The fact that $\alpha \in \Omega_f$ does follow from the discussion in [19] (using an identity from [31, p. 358]), but the result of Schertz [36] gives this more directly.  We will not make direct use of Proposition 2.2, but the relations in (2.10) are important for the arguments in Section 3, and curves $E_\alpha$ with complex multiplication will come up in the proof of Theorem 3.3.  \medskip

We will now show that $i \neq 0$ in (2.4). \bigskip

\noindent {\bf Lemma 2.3.} a) {\it If $w=k+\sqrt{-d}/2$ ($d$ even) or $(k+\sqrt{-d})/2$ ($d$ odd) is chosen so that
$$k^2 \equiv -d/4 \ \ \textrm{or}  \ -d \ (\textrm{mod} \ 9),\ \ k \equiv 1\ (\textrm{mod}\ 6),\eqno{(2.11)}$$ 
then the 24-th power $(\eta(w/9)/\eta(w))^{24}$ is a generator of the principal ideal $\wp_3'^{24}$ in the field $\Omega_f$, where $(3)=\wp_3 \wp_3'$ and $\wp_3=(3,w)$ in $R_K$.} \smallskip

\noindent b) {\it The principal ideal $(\alpha-3)=\wp_3'^3$, where $\wp_3'$ is the conjugate ideal of $\wp_3$ in $R_K$.  Furthermore, $(\alpha, 3)=\wp_3'$ and $\wp_3' || \alpha$ in $R_K$.}  \medskip

\noindent {\it Proof.} Part a) is classical and well-known.  We have that
\begin{eqnarray*}
\left(\frac{\eta(w/9)}{\eta(w)}\right)^{24}&=&\left(\frac{\eta(w/9)}{\eta(w/3)}\right)^{24}\left(\frac{\eta(w/3)}{\eta(w)}\right)^{24}\\
&=& \left(3^{12}\frac{\Delta(w,9)}{\Delta(w,3)}\right) \left(3^{12} \frac{\Delta(w,3)}{\Delta(w,1)}\right)= \varphi_{P_0}(w, 3) \varphi_{P_0}(w,1),
\end{eqnarray*}
in Hasse's notation [21, pp. 10-11], where 
$$\Delta(w_1,w_2)=\left(\frac{2\pi}{w_2}\right)^{12} \eta\left(\frac{w_1}{w_2}\right)^{24},$$
$$\varphi_{M}(w_1,w_2)=m^{12} \frac{\Delta(M(w_1,w_2))}{\Delta(w_1,w_2)}, \ \ \textrm{det}(M)=m,$$\smallskip
$M$ is a primitive, integral $2 \times 2$ matrix, and $P_0$ is the $2 \times 2$ diagonal matrix with diagonal entries 1 and 3. By the choice of $k$ we know that $\{w,9\}$ is a basis of the ideal $\wp_{3,-d}^2$.  Now Hasse's Satz 10 in [21] implies that both $\varphi_{P_0}(w, 3)$ and $\varphi_{P_0}(w, 1)$ are generators of the conjugate ideal $\wp_3'^{12}$, and this implies the lemma.  (See also [15, p.32] or [27, p. 165].) $\square$ \bigskip

We summarize the above discussion as follows. \bigskip

\noindent {\bf Theorem 2.4.} {\it If $-d=d_Kf^2 \equiv 1$ (mod $3$) and $w$ is given by
$$w=\cases{k+\frac{\sqrt{-d}}{2}, &if $2 \mid d$,\cr \frac{k+\sqrt{-d}}{2}, &if $2 \nmid d$,}$$
where $k$ satisfies (2.11), then for a unique value of $i \in \{1, 2\}$ the numbers 
$$\alpha = 3+\left(\frac{\eta(w/9)}{\eta(w)}\right)^3, \ \ \beta = 3\omega^i+ \omega^i \left(\frac{3\eta(3w)}{\eta(w/3)}\right)^3,$$
give a point $(\alpha,\beta)$ on the Fermat curve $Fer_3$ which is defined over the ring class field $\Omega_f$.} \medskip

\noindent {\it Proof.} To prove that $i \neq 0$ in the formula for $\beta$, note from [21, Satz 10] and [15, p. 43] that
$$\left(\frac{3\eta(3w)}{\eta(w/3)}\right)^{24}=\varphi_{P_4}(w,1) \frac{3^{12}}{\varphi_{P_0}(w,1)} \cong 1 \cdot \frac{3^{12}}{\wp_3'^{12}}=\wp_3^{12},$$
using the same notation as in Lemma 2.3, where $P_4$ is the diagonal matrix with entries $3$ and $1$ (and $\cong$ denotes divisor equality).  Hence, we have that
$$\left(\frac{3\eta(3w)}{\eta(w/3)}\right)^6 \cong \wp_3^{3} \ \ \textrm{in} \ \Omega_f(\omega).$$
This shows that $i = 0$ is impossible, since this would imply that $\left(\frac{3\eta(3w)}{\eta(w/3)}\right)^3 \in \Omega_f$, a domain in which $\wp_3^3$ is not the square of an ideal.  $\square$  \bigskip

\section{Properties of the solution.}

The main goal of this section is to complete the proof of Theorem 1 by showing that the number $\gamma =\omega^i\frac{\eta(w/9)}{\eta(w)}$ lies in $\Omega_f$ for some $i$.  This is needed for the proof of Theorem 4.2.  We begin with the following proposition. \bigskip

\noindent {\bf Proposition 3.1.} {\it a) The algebraic number $\alpha$ defined by $(2.4)$ has degree $2h(-d)$ over $\mathbb{Q}$ and therefore $\Omega_f=\mathbb{Q}(\alpha)$. \smallskip

b) The minimal polynomial $p_d(x)$ of $\alpha$ over $\mathbb{Q}$ is a normal polynomial. \smallskip

c) The linear fractional mapping $\sigma_1(x)=3(x+6)/(x-3)$ acts on the roots of $p_d(x)$. \smallskip

d) $Norm_{\Omega_f/\mathbb{Q}}(\alpha-3)=3^{3h(-d)}$.}\medskip

\noindent {\it Proof.} The degree of $\alpha$ over $\mathbb{Q}$ is divisible by $h=h(-d)$ since it generates $\Omega_f$ over $K$.  If $\alpha$ had degree $h$, then it would follow from $\mathbb{Q}(j_\alpha) \subseteq \mathbb{Q}(\alpha)$ and the fact that $j_\alpha$ has degree $h$ over $\mathbb{Q}$ that $\mathbb{Q}(j_\alpha) = \mathbb{Q}(\alpha)$; and since $j_\alpha=j(w)$ is real, $\alpha$ would be real.  Equation (2.10) would then imply that $j(w/3)$ is also real.  But $j(w/3)$ is the $j$-invariant of the ideal class of the prime divisor $\wp_{3,-d}$ of 3 in $\textsf{R}_{-d}$, and $\wp_{3,-d}$ has order $\le 2$ in the class group only if $9=N(\lambda)$ for some primitive element $\lambda=(a+b\sqrt{-d})/2$.  This can only be the case if $36 \ge d$, so $j(w/3)$ is not real for $d>36$.  The same conclusion follows if $d=23$, since in that case, $36=x^2+23y^2$ does not have a primitive solution.  Thus, $\alpha$ cannot have degree $h(-d)$, so it must have degree $2h(-d)$ since it lies in the field $\Omega_f$.  Otherwise $d < 36$ and $-d=-8,-11,-20, -32$, or $-35$.  In the first two cases $h(-d)=1$, and $j=8000$ and $-32768$, respectively.  Factoring the polynomial
$$F_d(x)=(x^3-27)^{h(-d)}H_{-d}\left(\frac{x^3(x^3-24)^3}{x^3-27}\right)$$
in these two cases shows that $\alpha$ is the root of an irreducible quadratic, verifying a).  If $d = 20, 32$ or $35$, then $h(-d)=2$, and factoring $F_d(x)$ shows that $\alpha$ is the root of an irreducible quartic.  The quadratics and quartics satisfied by $\alpha$ for these $5$ values of $d$ are given in equation (4.14) below.  (See [30] for the polynomials $H_{-d}(x)$ for $d=20, 32, 35$.)  This proves a) and b). \medskip

To prove c), note that $p_d(x)$ is irreducible over $\mathbb{Q}(\omega)$ since $\Omega_f \cap \mathbb{Q}(\omega)=\mathbb{Q}$.  From equation (2.10) it is clear that $p_d(x)$ is a factor of the polynomial
$$G_d(x)=(x^3-27)^{3h(-d)}H_{-d}\left(\frac{x^3(x^3+216)^3}{(x^3-27)^3}\right).\eqno{(3.1)}$$ \smallskip
Since the rational function
$$r(x)=\frac{x^3(x^3+216)^3}{(x^3-27)^3}$$
is mapped to itself under any transformation $\sigma=\left(x \rightarrow (ax+b)/(cx+d)\right)$ in the group $G_{12}=\textrm{Gal}(k(x)/k(r(x)))$, where $k=\mathbb{Q}(\omega)$, it follows that 
$$p_d^\sigma(x)=(cx+d)^{2h(-d)}p_d(\sigma(x))$$ \smallskip
is an irreducible factor of $G_d(x)$ over $\mathbb{Q}(\omega)$.  This is true for any irreducible factor of $G_d(x)$, so $G_{12}$ acts on irreducible factors of $G_d(x)$ over $\mathbb{Q}(\omega)$. \medskip

We show that $G_{12}$ is transitive on these factors (up to multiplication by nonzero constants).  This follows easily from the Galois theory of the normal extension $k(x)/k(r(x))$, but here is a direct proof.  For every root $j$ of $H_{-d}(x)$ there is a root $\xi$ of $p_d(x)$ for which $j=r(\xi)$, because $j(w/3)=r(\alpha)$ by (2.10), and there is an automorphism $\tau$ of $\Omega_f/K$ for which $j=j(w/3)^\tau$.  Thus, we can take $\xi=\alpha^\tau$, since the set of roots of $p_d(x)$ is mapped into itself by $\textrm{Gal}(\Omega_f/K)$.  If $q(x)$ is any irreducible factor of $G_d(x)$ over $\mathbb{Q}(\omega)$, and $\xi'$ is one of its roots, then $r(\xi')=j=r(\xi)$ is a root of $H_{-d}(x)$, so $\xi'$ and $\xi$ must belong to the same orbit under $G_{12}$ (since the $12$-th degree rational function $r(z)-r(w)$ factors into linear factors over $k$, where each factor of the numerator corresponds to an element of $G_{12}$).  From this it follows that $q(x)=c \cdot p_d^\sigma(x)$ for some element $\sigma \in G_{12}$ and some nonzero constant $c$. \medskip

Now, $\textrm{deg}(G_d(x)) = 12h$ and $\textrm{deg}(p_d(x))=2h$, so the stabilizer of $p_d(x)$ has order 2 and there is an involution in $G_{12}$ which fixes $p_d(x)$.  There are only three involutions in $G_{12}$, namely
$$\sigma_1(z)=\frac{3(z+6)}{z-3}, \hspace{.1 in} \sigma_3(z)=\frac{3\omega(z+6\omega)}{z-3\omega}, \hspace{.1 in} \sigma_4(z)=\frac{3\omega^2(z+6\omega^2)}{z-3\omega^2}.\eqno{(3.2)}$$
Suppose, for example, that $\sigma_3(\alpha)=s$ is a root of $p_d(x)$ and therefore lies in $\Omega_f$.  Then the cube root of unity $\omega$ satisfies
$$18\omega^2+(3\alpha+ 3s)\omega-s\alpha=0$$ \smallskip
over $\Omega_f$.  But $\omega^2 + \omega + 1=0$ implies that $18=3(\alpha+s)=-s\alpha$; hence, $\alpha$ is a root of $x^2-6x-18=0$ and $\alpha=3 \pm 3\sqrt{3}$, impossible since $\sqrt{3} \notin \Omega_f$.  The same argument shows that $\sigma_4(\alpha)$ cannot lie in $\Omega_f$.  Hence, $p_d(x)$ can only be fixed by $\sigma_1(z)$, and this proves c).  It follows from c) that the roots of $p_d(x)$ come in pairs $\xi, \sigma_1(\xi)=3+27/(\xi-3)$.  Therefore $\xi-3$ and $27/(\xi-3)$ are conjugates, so the constant term of $p_d(x+3)$ is $3^{3h(-d)}$, which proves d).  This fact is consistent with Lemma 2.3b).  $\square$ \bigskip

\noindent {\bf Corollary.} {\it The roots of the polynomial $G_d(x)$ are all contained in the field $\Omega_f(\omega)$.} \bigskip

Next, we prove the following relationship between $\alpha$ and $\beta$. \bigskip

\noindent {\bf Proposition 3.2.} {\it If $\tau = \left(\Omega_f/K ,\wp_3\right)$ is the Artin symbol for the prime ideal $\wp_3=(3,w)$ in $\Omega_f/K$, then
$$\alpha^\tau=\sigma_1(\beta)=\frac{3(\beta+6)}{\beta-3}$$ \smallskip
and $\beta$ is conjugate to $\alpha$ over $\mathbb{Q}$.}  \medskip

\noindent {\bf Remark.}  The Shimura Reciprocity Law, as presented in [36, p. 123], does not apply in this situation, since $\alpha=\mathfrak{f}(w)$, where $\mathfrak{f}$ is a modular function for $\Gamma(9)$, and $\tau$ is the Frobenius automorphism for $\wp_3$, which does not satisfy the condition of being relatively prime to the level ($9$) of $\mathfrak{f}$.  The formula of the theorem is equivalent to $\sigma_1(\mathfrak{f}(w))^\tau=\mathfrak{g}(w/3)^\tau=\omega^i \mathfrak{g}(w)$, for $i=1$ or $2$. \medskip

\noindent {\it Proof.} First note that the automorphism $\tau_1=(j(w/3) \rightarrow j(w))$ is equal to the Artin symbol $\tau=\left(\Omega_f/K ,\wp_3\right)$.  This follows from Hasse's congruence [15, p. 34], [21, Satz 11]: since $\{w, 1\}$ is a basis for $\textsf{R}_{-d}$, we have
$$j(w/3)^{\tau_1} = j(w) = j(\wp_{3,-d}^{-1} \wp_{3,-d})\equiv j(\wp_{3,-d})^3 =  j(w/3)^3\equiv j(w/3)^\tau \ \ (\textrm{mod} \hspace{.05 in} \wp_3),$$ \smallskip
which implies that $j(w/3)^{\tau_1} \equiv j(w/3)^\tau$ (mod $\wp_3$) and therefore $\tau_1=\tau$.  This uses the fact that 3 does not divide the discriminant of $H_{-d}(x)$.  (See [14].)  Secondly, Deuring's reduction theory [12] implies that 0 is not a root of $H_{-d}(x)$ (mod 3).  This is because $j=0$ is supersingular for $p=3$, but $\left(-d/3 \right)=+1$, so roots of $H_{-d}(x)$ (mod 3) are singular but not supersingular.  Hence 3 does not divide $H_{-d}(0)$, and the roots $j(w), j(w/3)$ of $H_{-d}(x)$ are relatively prime to 3. \medskip

It follows from this that $(\alpha, \beta, 3)=1$.  Suppose some prime divisor $\mathfrak{p}$ of 3 in $\Omega_f$ did divide $(\alpha, \beta)$.  From the formulas
$$j(w)=\frac{\alpha^3(\alpha^3-24)^3}{\alpha^3-27} \ \textrm{and} \ \frac{\alpha^3}{\alpha^3-27} = \frac{\beta^3}{27}$$
we have $\beta(\alpha^3-24)=3\xi$, where $\xi^3=j(w)$, so that $\mathfrak{p}^2$ would divide 3, impossible since 3 is not ramified in $\Omega_f/\mathbb{Q}$.  On the other hand, every prime divisor of 3 divides $\alpha$ or $\beta$, by the relation 
$$27\alpha^3+27\beta^3=\alpha^3\beta^3.\eqno{(3.3)}$$ \smallskip

From equation (2.10) we have that $j(w/3)=r(\alpha)$, and (3.3) gives that 
$$j(w)=\frac{\alpha^3(\alpha^3-24)^3}{\alpha^3-27}=\frac{\beta^3(\beta^3+216)^3}{(\beta^3-27)^3} = r(\beta).\eqno{(3.4)}$$
Therefore, $j(w)=j(w/3)^\tau = r(\alpha^\tau)=r(\beta)$ implies that $\alpha^\tau=\sigma(\beta)$ for some $\sigma \in G_{12}$.  It is easy to check that $\sigma(\beta) \notin \Omega_f$ unless $\sigma = 1$ or $\sigma_1$, because all other elements of $G_{12}$ have coefficients that involve $\omega$.  This is obvious for the substitutions
$$x \rightarrow \omega x, \hspace{.05 in} \frac{3 \omega (x+6)}{x-3}, \hspace{.05 in}\frac{3(x+6\omega)}{x-3\omega}$$ \smallskip
and the substitutions obtained from these by replacing $\omega$ by $\omega^2$.  For the maps $\sigma_3(x) = \frac{3\omega(x+6\omega)}{x-3\omega}$ and $\sigma_4(x)= \frac{3\omega^2(x+6\omega)}{x-3\omega}$, the claim follows by writing them as linear fractional expressions in $\omega$.  For example,
$$\frac{3\omega(x+6\omega)}{x-3\omega} =\frac{(3x-18)\omega-18}{-3\omega+x}.$$ \smallskip
The determinant, as a linear fractional expression in $\omega$, is $3(x^2-6x-18) \neq 0$ for $x \in \Omega_f$ (see the proof of Proposition 3.1), so $\sigma_3(\beta) \notin \Omega_f$.  (The argument for $\sigma_4(\beta)$ is similar, except the corresponding determinant is $-3(x+6)(x-3)$, and $\beta \neq 3,-6$, since $\alpha^3 \neq 24$.)  Thus, $\alpha^\tau=\sigma(\beta)$ for $\sigma =1$ or $\sigma_1$.\medskip

From these facts we can rule out the equation $\alpha^\tau=\beta$.  Suppose $\mathfrak{p}$ is any prime divisor of $\wp_3$ in $\Omega_f$.  Then $\mathfrak{p}$ is fixed by the automorphism $\tau=\left(\Omega_f/K, \wp_3\right)$, so that $\alpha^\tau=\beta$ would imply that $\mathfrak{p}$ divides $\alpha$ and $\beta$, contradicting $(\alpha, \beta, 3)=1$.  This proves that $\alpha^\tau=\sigma_1(\beta)$.  The fact that $\alpha$ and $\beta$ are conjugates follows from $\beta=\sigma_1(\alpha^\tau)$ and Proposition 3.1c). $\square$  \medskip

\noindent {\bf Corollary}.  {\it The polynomial $p_d(x)$ is the only irreducible factor of the polynomial $G_d(x)$ (see (3.1)) having roots in $\Omega_f$.} \medskip

\noindent {\bf Proof}. This follows immediately from Proposition 3.1c) and the statement that $\sigma(\beta) \notin \Omega_f$ unless $\sigma = 1$ or $\sigma_1$. \bigskip

We will now show that the solution $(\alpha, \beta)$ of $Fer_3$ in the ring class field $\Omega_f$ of $K$ has the property that $\alpha=3+\gamma^3$ and $\beta=3+\gamma'^3$ for certain elements $\gamma, \gamma' \in \Omega_f$. 
Let $\alpha'$ be any algebraic number for which the Deuring normal form $E_{\alpha'}$ has complex multiplication by $\textsf{R}_{-d}$.  If $\beta'$ is chosen so that $(\alpha', \beta')$ is a point on $Fer_3$, then $\beta'$ is a root of the polynomial $G_d(x)$ in (3.1), by (3.4) with $\alpha, \beta$ replaced by $\alpha', \beta'$.  By the corollary to Proposition 3.1, all the roots of $G_d(x)$ are contained in the field $\Omega_f(\omega)$.  Furthermore, for each root $j$ of $H_{-d}(x)=0$, there is an $\alpha'$ satisfying $j(E_{\alpha'})=j$, for which $\alpha', \beta' =\sigma_1(\alpha'^{\tau}) \in \Omega_f$, which we can see by just conjugating the solution of (2.9) that we constructed before by some element of $\textrm{Gal}(\Omega_f/K)$.  \medskip

In the following theorem we temporarily free $\alpha$ of its meaning in (2.4). \bigskip

\noindent {\bf Theorem 3.3.}  {\it Let $K=\mathbb{Q}(\sqrt{-d})$, with $-d \equiv 1$ (mod 3).  If $\alpha$ is any algebraic number for which $E_\alpha: Y^2+\alpha XY+Y = X^3$ has complex multiplication by the order $\textsf{R}_{-d}$, then $\alpha$ lies in the class field $\Sigma_9\Omega_f$ over $K$, where $\Sigma_9$ is the ray class field of conductor $9$ over $K$.  In particular, $\alpha$ generates an abelian extension of $K$.}  \medskip

\noindent {\it Proof.} Suppose that $\alpha$ is an arbitrary root of 
$$X^3(X^3-24)^3-j(X^3-27)=0,$$
for some root $j$ of $H_{-d}(x)=0$.  We choose another root $\alpha'$ of the same equation which lies in $\Omega_f$, and a $\beta' \in \Omega_f(\omega)$ for which $(\alpha', \beta')$ lies on $Fer_3$.  If $\alpha'=\omega^i \alpha$ for some $i$, then certainly $\alpha \in \Omega_f(\omega) = \Sigma_3\Omega_f \subseteq \Sigma_9\Omega_f$. \medskip

Otherwise, we will construct an explicit isomorphism $E_\alpha \cong E_{\alpha'}$.  By [29, Thm. 3.3] and the arguments of [29, pp. 262-263] we may take $\beta'=\sigma_1(\beta)=3(\beta+6)/(\beta-3) \in \Omega_f(\omega)$, by replacing $\beta$ by $\omega^i \beta$ for some $i$.  Then an isomorphism $\iota: E_\alpha \rightarrow E_{\alpha'}$ is given by the equations (3.17) in [29, Prop. 3.10]:
$$X'= -\frac{\gamma'}{\gamma}X+\gamma', \ \ Y'=\frac{\sqrt{-3}(\beta-3)}{9} \left(Y-\sqrt{-3} \omega^2 \frac{\delta}{\gamma} X-\omega \delta \right);\eqno{(3.5)}$$
where
$$\gamma = \frac{-3\beta}{\alpha(\beta-3)}, \quad \delta=\frac{\beta-3\omega}{\beta-3}$$ \smallskip
are the coordinates of a point $(X, Y)=(\gamma, \delta)$ of order 3 on $E_\alpha$, and
$$\gamma'=-\frac{\beta+6}{3\alpha'}$$ \smallskip
is the $X'$-coordinate of a point of order 3 on $E_{\alpha'}$.  Now, by the choice of $\alpha'$ and the fact that $\beta$ is a root of $G_d(x)=0$, the element $\gamma'$ lies in $\Omega_f(\omega) \subset \Sigma_9\Omega_f$.  Furthermore, $E_{\alpha'}$ has complex multiplication by $\textsf{R}_{-d}$, $(9,f)=1$, and the Weierstrass normal form of $E_{\alpha'}$ is
$$Y^2=4X^3-\left(\frac{\alpha'^4-24\alpha'}{12}\right)X+\frac{\alpha'^6-36\alpha'^3+216}{216}, \ \ \Delta=\alpha'^3-27,$$
so the values $g_2, g_3, \Delta$ for this curve lie in $\Omega_f$.  Note that $g_2g_3 \neq 0$, since $j$ is a root of $H_{-d}(x)$ and $-d \neq -3, -4$.  Thus, we may use a theorem of Franz [16], according to which the ray class field $\Sigma_9$ of conductor $\mathfrak{f}=(9)$ over $K=\mathbb{Q}(\sqrt{-d})$ satisfies
$$\Sigma_9 \Omega_f=K(j(E_{\alpha'}),h(E_{\alpha'}[9])),$$
\noindent where $\displaystyle h(P)=\frac{g_2g_3}{\Delta}x(P)$ is a Weber function on points of $E_{\alpha'}$.  This shows, that if $P'=(X', Y')$ is a point of order $9$ on $E_{\alpha'}$, then $X' =\frac{\Delta h(P')}{g_2g_3} \in \Sigma_9\Omega_f$.  It follows from (3.5) that $-X/\gamma=\alpha X (\beta-3)/(3\beta) \in \Sigma_9\Omega_f$ for any point $(X, Y)$ in $E_\alpha[9]$, which gives that $\alpha X \in \Sigma_9\Omega_f$.  \medskip

We now use [29, Prop.3.6, Remark], according to which the roots $x$ of the cubic equation
$$x^3-(3+\alpha)x^2+\alpha x+1 = 0$$ \smallskip
\noindent are the $X$-coordinates of points of order $9$ on $E_\alpha$.  It follows that
$$x^2(x-3-\alpha) = x^3-(3+\alpha)x^2 = -\alpha x -1 \in \Sigma_9\Omega_f.\eqno{(3.6)}$$ \smallskip
\noindent On the other hand, $\alpha^3=27\beta^3/(\beta^3-27)=r \in \Sigma_9\Omega_f$, so multiplying the inclusion (3.6) by $\alpha^3$ gives that $\alpha x-\alpha(3+\alpha) \in \Sigma_9\Omega_f$, and therefore $\alpha^2+3\alpha \in \Sigma_9\Omega_f$.  Now form the expression:
$$(\alpha^2+3\alpha)^2 -9(\alpha^2+3\alpha)-6r= \alpha^4-27\alpha =(r-27)\alpha.$$ \smallskip
\noindent This gives that $(r-27)\alpha \in \Sigma_9\Omega_f$.  But $r-27=\alpha^3-27=\Delta$ is the discriminant of the curve $E_\alpha$, which is non-zero, so we get that $\alpha \in \Sigma_9\Omega_f$.  $\square$ \bigskip

\noindent {\bf Remark.} Calculations on Maple suggest that the following stronger statement holds. \medskip

\noindent {\it Conjecture.} Let $K=\mathbb{Q}(\sqrt{-d})$, with $-d =d_Kf^2 \equiv 1$ (mod $3$) and $d_K$ the discriminant of $K$.  If $E_\alpha: \ Y^2+\alpha XY+Y = X^3$ has complex multiplication by the order $\textsf{R}_{-d}$ of discriminant $-d$ in $K$, then $\alpha$ lies in the ring class field $\Omega_9 \Omega_f=\Omega_{9f}$ of conductor $9f$ over $K$. \bigskip

Now we prove \bigskip

\noindent {\bf Theorem 3.4.}  {\it If $K=\mathbb{Q}(\sqrt{-d})$, with $-d =d_Kf^2 \equiv 1$ (mod 3), and $\alpha$ is defined by (2.4), then \medskip

\noindent a) $\alpha = 3 + \gamma^3$, where $\gamma \in \Omega_f$ and $(\gamma)=\wp_3'$. \smallskip

\noindent b) For some $i$, $\displaystyle \gamma = \omega^i \frac{\eta(w/9)}{\eta(w)}$, where $w$ is defined by (1.1).  \medskip

\noindent c) There is an integral point $(z,w)$ on the genus $19$ curve
$$C_{19}: \hspace{.1 in} z^3w^3(z^6+9z^3+27)(w^6+9w^3+27)=729 $$ \smallskip
with coordinates in the ring class field $\Omega_f$ of $K$.} \medskip

\noindent {\it Proof.} We use the same argument (suitably modified) as in the proof of [31, Theorem 1].  (Also see [28, Prop. 8.9].)  Let $k$ be field $k=\Sigma_9\Omega_f$, which by Theorem 3.3 contains the splitting field of the polynomial
$$F_d(x) = (x^3-27)^{h(-d)} H_{-d}\left( \frac{x^3(x^3-24)^3}{x^3-27} \right)$$ \smallskip
over $\mathbb{Q}$ (see (2.5)).  Also, $\Omega_f(\omega) \subset k$ and $k/\Omega_f$ is abelian.  \medskip

Now we let $\displaystyle j=j(\alpha)=\frac{\alpha^3(\alpha^3-24)^3}{\alpha^3-27}$, where $\alpha$ is an indeterminate over $k$.  We define $\beta$ by the condition $(\alpha, \beta) \in Fer_3$.  We will use the fact from [31] that the normal closure of the algebraic extension $k(\alpha)/k(j)$ is the function field $N$ given by
$$N=k(\alpha, \ (\beta-3)^{1/3}, \ (\omega \beta-3)^{1/3}, \ (\omega^2 \beta-3)^{1/3}).$$ \smallskip
\noindent Consider a root $j_0$ of the class equation $H_{-d}(x)$.  Then $j_0 \neq 0, 1728$ and all the roots of $j(\alpha)=j_0$ lie in the field $k$, by Theorem 3.3.  If $P_{j_0}$ is the prime divisor of the rational function field $k(j)$ corresponding to the polynomial $j-j_0$, this implies that $P_{j_0}$ splits into primes of degree 1 in the field $k(\alpha)$, and therefore splits completely in the normal closure $N$ of $k(\alpha)/k(j)$  ($j_0 \neq 0, 1728$ implies $P_{j_0}$ is unramified in $N/k(j)$).  Consider a root $\alpha_0$ of $F_d(x)$ for which $j(\alpha_0)=j_0$ and $\alpha_0, \beta_0 \in \Omega_f$, and any prime divisor $\mathfrak{P}$ of $N$ for which
$$\alpha \equiv \alpha_0 \quad \beta \equiv \beta_0 \ (\textrm{mod} \ \mathfrak{P}),$$ \smallskip
\noindent so that $\mathfrak{P}|P_{j_0}$.  Since $\mathfrak{P}$ has degree 1 over $k$, it follows that there is an element $\gamma_0 \in k$ for which
$$(\beta-3)^{1/3} \equiv \gamma_0  \ (\textrm{mod} \ \mathfrak{P}),$$ \smallskip
\noindent and therefore $\beta_0 \equiv \beta \equiv \gamma_0^3+3$ (mod $\mathfrak{P}$).  Hence, $\beta_0 = \gamma_0^3+3$ in $k$.  However, $\gamma_0$ generates an abelian extension of $\Omega_f$ and $\gamma_0^3=\beta_0-3 \in \Omega_f$.  This implies that $x^3-\gamma_0^3$ is reducible over $\Omega_f$: otherwise, its splitting field would have Galois group $S_3$ over $\Omega_f$, since $\omega \notin \Omega_f$.  Therefore, $\beta_0=3+\gamma_0^3$ for some $\gamma_0 \in \Omega_f$, and applying the automorphism $\tau^{-1} \in \textrm{Gal}(\Omega_f/K)$ to the equation $\sigma_1(\beta_0)=\alpha_0^\tau$ we get that
$$\alpha_0=\sigma_1(\beta_0^{\tau^{-1}}) = 3+\frac{27}{\beta_0^{\tau^{-1}}-3} = 3+\left(\frac{3}{\gamma_0^{\tau^{-1}}} \right)^3.$$ \smallskip
Therefore, $\alpha_0=3+\gamma^3$, with $\gamma=3/\gamma_0^{\tau^{-1}} \in \Omega_f$.  The remainder of part a) and part b) is immediate from this and Lemma 2.3b) and Theorem 2.4.  Finally, the curve $C_{19}$ arises from setting $\alpha=3+z^3, \beta=3+w^3$ in $Fer_3$ and simplifying.  Thus, $(z,w)=(\gamma,\gamma_0)$ is a point on $C_{19}$ defined over $\Omega_f$.  (See [31, Thm. 5].)  $\square$  \bigskip

This theorem verifies several assertions made in [31, p. 341].  In particular, part a) is an analogue of [31, Thm. 1] in characteristic $0$.  \medskip 

\noindent {\bf Remark.} Part b) of this theorem is related to a theorem of Fricke-Hasse-Deuring (see [17, III, p. 362], [22], or [15, p. 41]), which says that
$$\varphi_S(w_1,w_2)=(\textrm{det} \ S)^{12} \frac{\Delta(S(w_1,w_2)^t)}{\Delta(w_1,w_2)}$$
is the 24-th power of an element of $\Omega_f$.  Here $(w_1,w_2)$ is a basis of a proper ideal $\mathfrak{a}$ in $\textsf{R}_{-d}$, $S$ is a $2 \times 2$ integral matrix for which $S(w_1,w_2)^t$ is a basis of $\mathfrak{a} \mathfrak{b}^2$ with $\textrm{det} \ S >0$, and $(\mathfrak{b},6f)=1$.  In our case, $\mathfrak{b}=\wp_{3,-d}=\wp_3 \cap \textsf{R}_{-d}$, which is not relatively prime to $6$.  To handle this situation, in which $w/9$ is the ideal basis quotient for the ideal $\wp_{3,-d}^2$, one has to use the extension of this theorem due to Kubert and Lang [26, p. 296], and this extension only applies in the case $\Omega_1=\Sigma$.  The above proof avoids this complication.  \bigskip

Propositions 3.1 and 3.2 and Theorem 3.4 are summarized in Theorem 1 of the introduction, which is now completely proved. \bigskip

\noindent {\bf Proposition 3.5.} {\it Let $p_d(x)$ be the minimal polynomial of $\alpha$ over $\mathbb{Q}$, and $q_d(x)$ the minimal polynomial of the element $\gamma$ of $\Omega_f$ for which $\alpha=3+\gamma^3$.  Then we have the identities}
$$ (x-3)^{2h(-d)} p_d\left(\frac{3(x+6)}{x-3} \right)=3^{3h(-d)}p_d(x),$$
$$x^{2h(-d)} q_d \left(\frac{3}{x} \right) = 3^{h(-d)} q_d(x).\eqno{(3.7)}$$ \medskip

\noindent {\it Proof.} We know the mapping $\displaystyle x \rightarrow \sigma_1(x)=\frac{3(x+6)}{x-3} = 3+\frac{27}{x-3}$ permutes the roots of $p_d(x)$, so 
$$ (x-3)^{2h(-d)} p_d\left(\frac{3(x+6)}{x-3} \right)=c \cdot p_d(x),$$
for some constant $c$.  Now put $x=3+3\sqrt{3}$, a fixed point of the linear fractional map $\sigma_1(x)$.  This gives
$$(3\sqrt{3})^{2h(-d)} p_d(3+3\sqrt{3})=c \cdot p_d(3+3\sqrt{3}).$$
The quantity $p_d(3+3\sqrt{3})$ cannot be 0, since $3$ does not ramify in $\Omega_f$.  Hence, $c=3^{3h(-d)}$.  We also know that $p_d(3+x^3)=q_d(x)q_d(\omega x)q_d(\omega^2 x)$.  Using $\displaystyle \sigma_1(3+x^3)=3+\frac{27}{x^3}$ gives that
$$x^{6h(-d)}p_d\left(3+\frac{27}{x^3} \right)=3^{3h(-d)}p_d(3+x^3),$$
hence that

$$x^{6h(-d)}q_d\left(\frac{3}{x} \right)q_d\left(\frac{3\omega}{x} \right)q_d\left(\frac{3\omega^2}{x} \right) = 3^{3h(-d)}q_d(x)q_d(\omega x)q_d(\omega^2 x).$$ \smallskip

\noindent Since the roots of $x^{2h(-d)}q_d(3/x)$ lie in $\Omega_f$, we must have $x^{2h(-d)}q_d\left(\frac{3}{x} \right)=c_1 q_d(x)$, with some constant $c_1$, and putting $x=\sqrt{3}$ shows immediately that $c_1=3^{h(-d)}$.  $\square$

\section{Solutions of $Fer_3$ as periodic points.}

\noindent The substitution
$$x=\frac{9\beta}{\alpha(\beta-3)}, \ \ y=\frac{9\beta}{\beta-3}\eqno{(4.0)}$$
converts the solution $(\alpha,\beta)$ of $Fer_3$ into a point $(x,y)$ on the elliptic curve
$$E: \hspace{.1 in} Y^2-9Y=X^3-27.\eqno{(4.1)}$$ \smallskip
(Putting $-Y+4$ for $Y$ in (4.1) shows this is the curve (27B) in [11].)  Proposition 3.2 shows that the coordinates of the point $(x,y)$ on the new curve are
$$x = \frac{\alpha^\tau + 6}{\alpha}, \hspace{.1 in} y = \alpha^\tau + 6=\alpha x.\eqno{(4.2)}$$ \smallskip
In particular, the point $(x,y)$ has integral coordinates.  Note also that
$$E(\mathbb{Q})=\{O, (3,0), (3,9) \},$$
and the points $(3,0)$ and $(3,9)$ have order $3$ on $E$.  \medskip

We write the equation expressing the fact that the point $P_d=(\frac{\alpha^\tau+6}{\alpha},\alpha^\tau+6)$ lies on $E$ in the following form:
$$(\alpha^\tau+6)^2-9(\alpha^\tau+6)-\frac{(\alpha^\tau+6)^3}{\alpha^3}+27=\frac{((\alpha^\tau)^2+3\alpha^\tau+9)\alpha^3-(\alpha^\tau+6)^3}{\alpha^3}=\frac{g(\alpha,\alpha^\tau)}{\alpha^3},$$ \smallskip
where
$$g(x,y)=(y^2+3y+9)x^3-(y+6)^3.\eqno{(4.3)}$$ \smallskip
Now define the following functions:
$$T(z)=\frac{z^2}{3}(z^3-27)^{1/3}+\frac{z}{3}(z^3-27)^{2/3}+\frac{z^3}{3}-6,\eqno{(4.4)}$$
$$S(z)=\frac{z+6}{(z^2+3z+9)^{1/3}}.\eqno{(4.5)}$$ \smallskip
These functions satisfy the equations
$$g(x,T(x))=0, \quad g(S(y),y)=0$$ \smallskip
so they are solutions of the implicit relation $g(x,y)=0$ and inverse functions of each other.  Note that the formula for $T(z)$ can be obtained by using the Ferrarro-Tartaglia-Cardan formulas to solve $g(z,y)=0$ as a cubic in $y$.  \medskip

We will take $T(z)$ and $S(z)$ to be defined on certain subsets of the maximal unramified extension $\textsf{K}_3$ of $\mathbb{Q}_3$ inside its algebraic closure $\bar \mathbb{Q}_3$.  This is permissible because the cube root $(z^3-27)^{1/3}$ can be defined as the convergent $3$-adic series

$$(z^3-27)^{1/3}=z \sum_{k=0}^\infty{{\frac{1}{3} \atopwithdelims ( ) k} \left( \frac{-3}{z} \right)^{3k}}, \quad w_3(z) \le 0,$$ \smallskip

\noindent for $z$'s satisfying $3w_3(-3/z)>3/2$, i.e. whenever $|z|_3 \ge 1$ in  $\textsf{K}_3$. (See [34, p. 34].)  Thus we have

$$T(z)=-6+z^3 + \frac{z^3}{3} \sum_{k=1}^\infty{ \left({\frac{1}{3} \atopwithdelims ( ) k}+{\frac{2}{3} \atopwithdelims ( ) k} \right) \left( \frac{-3}{z} \right)^{3k}}$$
$$=z^3-15-z^3 \sum_{k=2}^\infty{  \frac{3^{2k-1} b_k}{k!} \frac{1}{z^{3k} }}, \ \ |z|_3 \ge 1,\eqno{(4.6)}$$ \smallskip

\noindent where $b_k$ are the rational integers defined by $b_1=3$ and
$$b_k=(3k-4)(3k-7) \cdots 5 \cdot 2 + 2 \cdot (3k-5)(3k-8) \cdots 4 \cdot 1, \quad k \ge 2.$$ \smallskip
\noindent We note that $3^{k-1}b_k/k! \in \mathbb{Z}$ for $k \ge 1$ since $ {1/3 \atopwithdelims ( ) k}$ and $ {2/3 \atopwithdelims ( ) k}$  are $p$-adic integers for any prime $p \neq 3$ and since $w_3(3^k/k!) > w_3(3)=1$ implies that $w_3(3^{k-1}/k!)  > 0$.  (See [24, pp. 264-265].)  The field $\textsf{K}_3$ is not complete with respect to its valuation, but every finite subextension is complete, so the series for $T(z)$ converges to an element of the field $\mathbb{Q}_3(z)$.  \medskip

It is clear from (4.6) that $T(z)$ is a $3$-adic unit whenever $z$ is, and that 
$$T(z) \equiv -6+z^3 \quad (\textrm{mod} \ 9), \quad w_3(z)=0.\eqno{(4.7)}$$  \smallskip
More generally, (4.6) implies that
$$|T(z)-z^3|_3<1, \ \textrm{for} \ |z|_3 \ge 1,$$
and therefore
$$|T(z)|_3=|z|_3^3, \ \textrm{for} \ |z|_3 > 1.\eqno{(4.8)}$$
Hence, we may iterate the function $T(z)$ on the set $\{z: |z|_3 \ge 1\} \subset \textsf{K}_3$.  In particular, (4.8) shows that $T(z)$ has no periodic points in the region $|z|_3 >1$, so all of its periodic points must lie in the unit group $\textsf{U}_3$ of $\textsf{K}_3$.  \medskip

On the other hand, the function $S(z)$ is convergent in a disc about $z=3$, since
$$S(3+27z)=\frac{3(1+3z)}{(1+9z+27z^2)^{1/3}} =3(1+3z) \sum_{k=0}^\infty{{-\frac{1}{3} \atopwithdelims ( ) k} 3^{2k} z^k (1+3z)^k}$$
$$=3+27z^3-243z^4+1458z^5-6804z^6+\cdots \equiv 3 \quad (\textrm{mod} \hspace{.05 in} 27), \quad w_3(z) \ge 0.$$\smallskip
This shows that we may iterate $S(z)$ on the disc $\textsf{D}_3=\{z: |z-3|_3 \le |3|_3^{3} \}$. \bigskip

\noindent {\bf Lemma 4.1.} {\it If $z \in \textsf{K}_3$ satisfies $|z|_3 \ge 1$ and $g(z,w)=0$ for some $w \in \textsf{K}_3$, then $w=T(z)$.  Similarly, $z \in \textsf{D}_3$ and $g(w,z)=0$ for $w \in \textsf{K}_3$ implies that $w=S(z)$.}  \medskip

\noindent {\it Proof.} Let $t=T(z)$. Then $g(z,t)=g(z,w)=0$ implies that 

$$0=\frac{(w+6)^3g(z,t)-(t+6)^3g(z,w)}{z^3}=(t^2+3t+9)(w+6)^3-(w^2+3w+9)(t+6)^3$$
$$=(w-t) \left( (t^2+3t+9)w^2+(3t^2-45t-54)w+9t^2-54t+324 \right).$$ \smallskip
If $w \neq t$, then $w$ would be a root of the quadratic in the last equation, which has discriminant $\delta=-27(t+6)^2(t-3)^2$.  This would imply that $\mathbb{Q}_3(w,t)=\mathbb{Q}_3(\sqrt{-3},t)$ is contained in $\textsf{K}_3$, which is impossible.  The second assertion follows in the same way. $\square$ \bigskip

Since $-d \equiv 1$ (mod 3), $-d$ has a square root in $\mathbb{Q}_3$, and any ring class field $\Omega_f$ of $K=\mathbb{Q}(\sqrt{-d})$ whose conductor is prime to 3 embeds into the field $\textsf{K}_3$.  There are two embeddings of $K$ in $\mathbb{Q}_3$ corresponding to the completions $K_{\wp_3}$ and $K_{\wp_3'}$.  The functions $T(z)$ and $S(z)$ allow us to express the relation $g(\alpha,\alpha^\tau)=0$ as $g(\alpha,T(\alpha))=0$ or $g(\alpha^{\tau^{-1}},\alpha)=g(S(\alpha),\alpha)=0$ depending on which embedding of $\Omega_f$ we are considering:
$$g(\alpha,\alpha^\tau)=0 \hspace{.05 in} \textrm{and} \hspace{.05 in} \alpha \in \textsf{U}_3 \Rightarrow \alpha^\tau=T(\alpha), \eqno{(4.9)}$$
$$g(\alpha,\alpha^\tau)=0 \hspace{.05 in} \textrm{and} \hspace{.05 in} \alpha \in \textsf{D}_3 \Rightarrow \alpha^{\tau^{-1}}=S(\alpha). \eqno{(4.10)}$$ \smallskip
Recall that if $\mathfrak{p}$ is a prime divisor of $\wp_3$ or $\wp_3'$ in $\Omega_f$, $L$ is the decomposition field of $\mathfrak{p}$, and $\mathfrak{p}_L$ is the prime divisor of $L$ which $\mathfrak{p}$ divides, then $\textrm{Gal}(\Omega_f/L) \cong \textrm{Gal}(\left(\Omega_f\right)_\mathfrak{p}/L_{\mathfrak{p}_L}) = \textrm{Gal}(\left(\Omega_f\right)_\mathfrak{p}/\mathbb{Q}_3)$ is generated by $\tau$.  Thus we may apply $\tau$ to elements of $\left(\Omega_f\right)_\mathfrak{p} =\mathbb{Q}_3(\alpha) \subset \textsf{K}_3$, and we have
$$T(z)^\tau=T(z^\tau), \ \textrm{for} \ z \in \mathbb{Q}_3(\alpha), \ |z|_3 \ge 1,$$
$$S(z)^{\tau^{-1}}=S(z^{\tau^{-1}}),  \ \textrm{for} \ z \in \mathbb{Q}_3(\alpha) \cap \textsf{D}_3.$$
Then $\tau^n=1$ implies that
$$\alpha=\alpha^{\tau^{n}} = T^n(\alpha),\quad (\alpha, \mathfrak{p})=1,$$ 
$$\alpha=\alpha^{\tau^{-n}} = S^n(\alpha),\quad \mathfrak{p}|\alpha.$$ \smallskip
Hence, the solutions of $Fer_3$ we have constructed in ring class fields over $K$ (with conductors prime to 3) are periodic points of the algebraic functions $S$ and $T$! \medskip

We may find the minimal polynomials $p_d(x)$ of the periodic points of $S$ and $T$ using iterated resultants.  To find the points of period $n$, we look for the integers $d\equiv 2$ (mod $3$) for which the Frobenius automorphism $\tau=(\Omega_f/K,\wp_3)$ has order $n$.  Then the point $(\alpha,\alpha^\tau)$ is a point on the curve $g(x,y)=0$.  Applying $\tau$ repeatedly to this point gives the equations 
$$g(\alpha,\alpha^\tau)=g(\alpha^\tau,\alpha^{\tau^2})= \cdots = g(\alpha^{\tau^{n-1}},\alpha)=0, \ \ \tau^n=1.\eqno{(4.11)}$$ \smallskip
Thus, $\alpha$ will be the root of a series of nested resultants. \medskip

As an example, we find the periodic points of period $3$.  The condition $\tau^3=1$ implies the equations
$$g(\alpha,\alpha^\tau)=g(\alpha^\tau,\alpha^{\tau^2})=g(\alpha^{\tau^2},\alpha)=0.$$ \smallskip
\noindent Now compute the double resultant
$$R_3(x)=Res_{x_2}(Res_{x_1}(g(x,x_1),g(x_1,x_2)),g(x_2,x))=(x-3)(x^2+4x+6)(x^2+2x+12)$$
$$\times (x^6+11x^5+65x^4+191x^3+441x^2+405x+675)$$
$$\times(x^6+20x^5+126x^4+172x^3+180x^2-1188x+1188)$$
$$\times(x^6+22x^5+208x^4-40x^3+144x^2-3456x+6912)$$
$$\times(x^6+6x^5+560x^4-1384x^3+576x^2-12960x+43200)$$
$$\times(x^6-74x^5+1680x^4-6184x^3+2736x^2-43200x+172800)$$
$$\times(x^6-13x^5+841x^4-2567x^3+1071x^2-20493x+75141)$$
$$\times(x^{12}-44x^{11}+724x^{10}+11008x^9+30440x^8-125456x^7-806960x^6-1971936x^5$$
$$+4056480x^4+17611776x^3+46267200x^2+10730880x+24681024).$$ \smallskip

By our theory, every ring class field $\Omega_f$ for which $-d \equiv 1$ (mod 3) and $\tau^3=1$ must show up as the splitting field of one of the factors of this double resultant.  Note that the factors $x^2+4x+6$ and $x^2+2x+12$ have the respective discriminants $-d=-8$ and $-2^2d=-44$ and the corresponding fields $K$ have class number 1, i.e. $\tau=1$.  Further, $z=3$ is a fixed point of the mapping $S(z)$.\medskip

The sextic factors of $R_3(x)$ are the minimal polynomials $p_d(x)$ of $\alpha$ corresponding to the values of $d = 23, 44, 59, 83, 107, 92$, respectively, and the 12th degree polynomial is the minimal polynomial of $\alpha$ for $d=104$.  The second sextic is the polynomial $p_{44}(x)$, whose splitting field is the ring class field with conductor $(2)$ over $\mathbb{Q}(\sqrt{-11})$.  The sixth polynomial $p_{92}(x)$ corresponds to the ring class field with conductor (2) over $K=\mathbb{Q}(\sqrt{-23})$, which in this case coincides with the Hilbert class field of $K$. (This gives a second set of integral solutions of $Fer_3$ in this field.) Thus there are exactly 4 quadratic fields $K$ with $-d \equiv 1$ (mod 3) and $h(K)=3$, namely:
$$h(K)=3 \ \textrm{and} \ \textrm{ord}(\tau)=3 \ \ \textrm{iff} \ \ d = 23, 59, 83, 107.\eqno{(4.12)}$$ \smallskip
Taking $n=1$, note that
$$R_1(x)=g(x,x)=(x-3)(x^2+4x+6)(x^2+2x+12),\eqno{(4.13)}$$ \smallskip
so no fields with $h(K)>1$ have $\tau=1$.  We also note that the above computation implies that there is only one field $K=\mathbb{Q}(\sqrt{-d})$ of the required form, namely $d=104$, for which $h(K)=6$ and the Frobenius automorphism $\displaystyle \tau=(\Sigma/K,\wp_3)$ has order 3.  \medskip

Similarly, the resultant
$$R_2(x)=Res_{x_1}(g(x,x_1),g(x_1,x)) = -(x-3)(x^2+4x+6)(x^2+2x+12)$$
$$\times (x^4-12x^3+28x^2+48x+576)(x^4-8x^3+26x^2+60x+450)$$
$$\times (x^4+2x^3+26x^2+60x+180) \eqno{(4.14)}$$ \smallskip
\noindent shows there are no fields $K$ with $h(K)=6$ for which $\tau$ has order 2.  So in all cases where $h(K)=6$ except the case $d=104$, the automorphism $\tau$ has order 6.  (The values of $d$ corresponding to the the three quartics in this factorization are, respectively, $d=35$, $d=32$, and $d=20$.)  \medskip

In the case of $n=3$, one could, of course, determine the values of $d$ in (4.12) by solving the equations $4^r \cdot 3^3 = x^2+dy^2$ for $r=0,1$, but computing the double resultant $R_3(x)$ immediately gives minimal polynomials for the generators of the corresponding class fields! \medskip

We apply these insights to generalize the factorizations of $R_n(x)$ for $n = 1, 2, 3$.  We define $R_n(x)$ as follows.  First define $R^{(1)}(x,x_1)=g(x,x_1)$ and then recursively define
$$R^{(k)}(x,x_k)=Res_{x_{k-1}}(R^{(k-1)}(x,x_{k-1}),g(x_{k-1},x_k)), \quad k \ge 3.\eqno{(4.15)}$$
Then we set $x_n=x$ in $R^{(n)}(x,x_n)$ to obtain $R_n(x)$:
$$R_n(x) = R^{(n)}(x,x), \ \ n \ge 1.\eqno{(4.16)}$$ \smallskip
From this definition it is easy to see that the roots of $R_n(x)$ are exactly the $a$'s for which there exist common solutions of the equations
$$g(a,a_1)=0, \quad g(a_1,a_2)=0, \quad \cdots \quad g(a_{n-1},a)=0.\eqno{(4.17)}$$ \smallskip
If $a \in \textsf{U}_3$ or $\textsf{D}_3$, then 
$$a=T(a_{n-1})=T^2(a_{n-2})=\cdots=T^{n-1}(a_1)=T^n(a)$$
or
$$a=S(a_1)=S^2(a_2)=\cdots = S^{n-1}(a_{n-1})=S^n(a),$$\smallskip
respectively.  The roots of $R_n(x)$ in $\textsf{U}_3 \cup \textsf{D}_3$ are therefore points in $\textsf{K}_3$ whose periods with respect to $S$ or $T$ divide $n$. We will now show that all the roots of $R_n(x)$ are in $\textsf{U}_3 \cup \textsf{D}_3$.  \medskip

We first compute the total number of roots of $R_n(x)$ (with multiplicities) by considering the polynomial $g(x,x_1) \equiv x_1^2x^3-x_1^3\equiv x_1^2(x^3-x_1)$ (mod 3).  The definition of the resultant gives modulo $3$ that

$$R^{(2)}(x,x_2)\equiv Res_{x_1}(x_1^2(x^3-x_1),x_2^2(x_1^3-x_2))\equiv -x_2^8(x^9-x_2),$$
$$R^{(3)}(x,x_3)\equiv Res_{x_2}(-x_2^8(x^9-x_2), x_3^2(x_2^3-x_3))\equiv x_3^{26}(x^{27}-x_3),$$
and recursively
$$R^{(i+1)}(x,x_{i+1})\equiv Res_{x_i}\left( (-1)^{i-1} x_i^{3^i-1}(x^{3^i}-x_i),x_{i+1}^2(x_i^3-x_{i+1}) \right)$$
$$\equiv (-1)^i x_{i+1}^{3^{i+1}-1}(x^{3^{i+1}}-x_{i+1}) \ \ (\textrm{mod} \ 3).$$ \smallskip
Hence, taking $i=n-1$ and $x_n =x$ we have
$$R_n(x)\equiv (-1)^{n-1} x_n^{3^n-1}(x^{3^n}-x_n) \equiv (-1)^{n-1} x^{3^n-1}(x^{3^n}-x) \quad (\textrm{mod} \hspace{.05 in} 3).\eqno{(4.18)}$$ \smallskip
\noindent We deduce easily from (4.15) and the fact that the highest degree terms in $g(x,y)$ in $x$ and $y$ do not vanish (mod $3$) that
$$\textrm{deg}(R_n(x))=2 \cdot 3^n - 1, \quad n \ge 1.$$ \smallskip
It is clear from (4.17) that the roots of $R_k(x)$ are all roots of $R_n(x)$ whenever $k \mid n$.  We will show that $R_n(x)$ has distinct roots for all $n \ge 1$.  It will follow that the points of primitive period $n$ are the roots of a polynomial $\textsf{P}_n(x)$, for which
$$R_n(x) = \pm \prod_{k \vert n}{\textsf{P}_k(x)},\eqno{(4.19a)}$$
and therefore
$$\textsf{P}_n(x)=\pm \prod_{k|n}{R_k(x)^{\mu(n/k)}}.\eqno{(4.19b)}$$ \medskip

From (4.18) and Hensel's Lemma it is clear that for each irreducible factor $\bar f_i(x)$ of degree $n$ of $R_n(x)$ (mod 3) (for $n \ge 2$) there exists a monic irreducible polynomial $f_i(x)$ of degree $n$ in $\mathbb{Z}_3[x]$, dividing $R_n(x)$, for which $f_i(x) \equiv \bar f_i(x)$ (mod 3).  Moreover, if $\zeta$ is a root of $f_i(x)$ in $\bar \mathbb{Q}_3$, then $\mathbb{Q}_3(\zeta)$ is an unramified extension of $\mathbb{Q}_3$ of degree $n$, and therefore $\zeta \in \textsf{K}_3$.  By (4.18), $\zeta$ is not a root of $R_k(x)$ with $k < n$, so the fact that $\zeta \in \textsf{U}_3$ (for $n>1$) implies that $\zeta$ is a periodic point of $T$ with primitive period $n$.  It follows that $T(\zeta)$ is also a periodic point of primitive period $n$, and from (4.7) we know that
$$T(\zeta) \equiv \zeta^3 \quad (\textrm{mod} \ 3).$$ \smallskip
\noindent It follows from this congruence that $T(\zeta)$ and $\zeta$ are roots of the same $\bar f_i(x)$ (mod 3), and therefore they must be roots of the same factor $f_i(x)$, since (4.18) shows that the roots of $\bar f_i(x)$ are simple roots of $R_n(x)$ (mod 3).  Hence, the map
$$\zeta \rightarrow T(\zeta), \quad f_i(\zeta)=0$$ \smallskip
represents an automorphism of the field $\mathbb{Q}_3(\zeta)/\mathbb{Q}_3$.  The function $T(z)$ is therefore a lift of the Frobenius automorphism to $\mathbb{Q}_3(\zeta)$ (when applied to the roots of $f_i(x)$).  This argument shows that $R_n(x)$ has the $N_3(n)$ simple factors $f_i(x)$ of degree $n$ whose roots are units in $\textsf{K}_3$, where $N_3(n)$ is the number of monic irreducible polynomials of degree $n$ in $\mathbb{F}_3[x]$. \medskip

Now we consider the numbers $\displaystyle \sigma_1(\zeta)=\frac{3(\zeta+6)}{\zeta-3}=3+\frac{27}{\zeta-3}$, for a root $\zeta$ of any of the polynomials $f_i(x)$.  Since $\zeta \in \textsf{U}_3$, it is clear that $\sigma_1(\zeta) \in \textsf{D}_3-\{3 \}$.  We claim that $\sigma_1(\zeta)$ is a periodic point of $S(z)$ with primitive period $n$, and therefore also a root of $R_n(x)$.  To show this we will use the identity
$$S(\sigma_1(\zeta))=\sigma_1(T(\zeta)).\eqno{(4.20)}$$ \smallskip
To prove this, let $z=\zeta$ and $t=T(\zeta)$ in the identity
$$(\sigma_1(z)+6)^3(t-3)^3-27(t+6)^3(\sigma_1(z)^2+3\sigma_1(z)+9)$$
$$=\frac{-3^9}{(z-3)^3} \left((t^2+3t+9)z^3-(t+6)^3 \right)=\frac{-3^9}{(z-3)^3} g(z,t).$$
This gives
$$S(\sigma_1(\zeta))^3=\sigma_1(T(\zeta))^3.$$ \smallskip
The cube roots of unity do not lie in $\textsf{K}_3$, so we may take cube roots in the last equation, yielding (4.20).  This identity implies that the maps $S$ and $T$ are conjugate maps on their respective domains.  It follows that
$$S^n(\sigma_1(\zeta))=S^{n-1}(\sigma_1(T(\zeta)))=S^{n-2}(\sigma_1(T^2(\zeta)))= \cdots =\sigma_1(T^n(\zeta))=\sigma_1(\zeta).$$ \smallskip
\noindent Hence $\sigma_1(\zeta)$ is a periodic point of $S$ with primitive period dividing $n$.  A similar argument shows that if $\sigma_1(\zeta)$ had period less than $n$, then the period of $\zeta$ with respect to $T$ would also be less than $n$.  Hence, $\sigma_1(\zeta)$ is a root of $R_n(x)$ but not a root of $R_k(x)$ for any proper divisor $k$ of $n$.  Since $\sigma_1(z)$ is a linear fractional map it is clear that $\sigma_1: \textsf{U}_3 \rightarrow \textsf{D}_3$ is 1-1.  Furthermore, the sets $\textsf{U}_3$ and $\textsf{D}_3$ are disjoint.  This yields an additional $nN_3(n)$ distinct primitive roots of $R_n(x)$, and therefore a total of $2nN_3(n)$ primitive roots.  Since roots of $R_k(x)$ are also roots of $R_n(x)$ for $k \mid n$, and $x-3$ is the only factor of $R_1(x)$ in (4.13) to which the above analysis does not apply, it follows that the number of distinct roots of $R_n(x)$ is at least
$$\sum_{k \vert n}{2kN_3(k)}-1=2\cdot 3^n-1=\textrm{deg}(R_n(x)).$$
Hence, all the roots of $R_n(x)$ are simple.  This allows us to define the polynomial $\textsf{P}_n(x)$ by (4.19), and the degree of $\textsf{P}_n(x)$ is given by
$$\textrm{deg}(\textsf{P}_n(x))=\sum_{k \vert n}{\mu(n/k) \ \textrm{deg}(R_k(x)) }=2nN_3(n), \quad n >1.\eqno{(4.21)}$$
By the above arguments we also know that $\textsf{P}_n(x)$ factors over $\mathbb{Q}_3$ as
$$\textsf{P}_n(x)=\prod_i{f_i(x) \tilde f_i(x)},\eqno{(4.22)}$$
where the roots of $\tilde f_i(x)$ are the images of the roots of $f_i(x)$ under $\sigma_1$.  Therefore, all irreducible factors of $\textsf{P}_n(x)$ over $\mathbb{Q}_3$ have degree $n$.  (By (4.13) this is even true for $\textsf{P}_1(x)=R_1(x)$.)  Moreover, we may write the polynomial $f_i(x)$ in the form
$$f_i(x)=\prod_{k=0}^{n-1}{(x-T^k(\zeta_i))},$$
where $\zeta_i$ is a representative of the orbit under $T$ whose elements are the roots of $f_i(x)$.  \medskip

Now it is obvious from the definition (4.19b) of the polynomial $\textsf{P}_n(x)$ that it has coefficients in $\mathbb{Q}$.  Furthermore, we also know from (4.9) that if $(\alpha,\beta)$ is the solution of $Fer_3$ we constructed in the ring class field of $K$ corresponding to the discriminant $-d$ and $p_d(x)$ is the minimal polynomial of $\alpha$ over $\mathbb{Q}$, then
$$\textrm{ord}(\tau)=n \quad \Rightarrow \quad p_d(x) \hspace{.02 in} | \hspace{.02 in} \textsf{P}_n(x).$$ \smallskip
\noindent Conversely, every root $\zeta$ of $\textsf{P}_n(x)=0$ is of course algebraic over $\mathbb{Q}$, and by (4.2) and (4.10) either $(\frac{T(\zeta)+6}{\zeta},T(\zeta)+6)$ or $(\frac{\zeta+6}{S(\zeta)},\zeta+6)$ is a point on the elliptic curve $E$, because $g(\zeta,T(\zeta))=0$ or $g(S(\zeta),\zeta)=0$.  Replacing $\zeta$ in the second of these points by $\sigma_1(\zeta)$ for a 3-adic unit $\zeta$, these points convert via the inverse of (4.0) to the points $(\alpha,\beta)=\left( \zeta,\sigma_1(T(\zeta)) \right)$ and $(S(\sigma_1(\zeta)),\zeta)=(\beta, \alpha)$ on $Fer_3$ in an extension $L=\mathbb{Q}(\alpha, \beta) \subset \textsf{K}_3$ which is unramified over the prime $3$.   \medskip

The factorization of the resultants $R_n$ for $1 \le n \le 5$ suggests the following theorem. \bigskip

\noindent {\bf Theorem 4.2.} {\it For $n>1$, the polynomial $\textsf{P}_n(x)$ in (4.22) is the product of the polynomials $p_d(x)$ over all discriminants $-d \equiv 1$ (mod 3) for which the Frobenius automorphism $\tau=(\Omega_f/K,\wp_3)$ has order $n$.  This fact is equivalent to the formula}
$$\sum_{-d \equiv 1 (3), \hspace{.01 in}  ord(\tau)=n}{h(-d)}=nN_3(n)=\sum_{k|n}{\mu (n/k)3^k}, \quad n>1.\eqno{(4.23)}$$ \medskip

\noindent {\it Proof.} Let $\zeta \in \textsf{U}_3$ be a root of $\textsf{P}_n(x)$ in $\textsf{K}_3$.  As stated above, the point  $(\alpha,\beta)=\left( \zeta,\sigma_1(T(\zeta)) \right)$ lies on $Fer_3$.  This implies that there is an isogeny $\phi: E_\alpha \rightarrow E_\beta$ of degree $3$, where $E_\alpha$ is the curve (2.1) in Deuring normal form.  Let $\tau$ denote the automorphism $\tau: \zeta \rightarrow T(\zeta)$ on $\mathbb{Q}_3(\zeta)$.  By (2.5) and the fact that $(\alpha,\beta)$ lies on $Fer_3$ we have
$$j(E_\beta)=\frac{\beta^3(\beta^3-24)^3}{\beta^3-27}=\frac{\alpha^3(\alpha^3+216)^3}{(\alpha^3-27)^3}.\eqno{(4.24)}$$ \smallskip
Since the last rational function in this equation is invariant under the substitution $\alpha \rightarrow \sigma_1(\alpha) = \sigma_1(\zeta)=\beta^{\tau^{-1}}$, we have
$$j(E_\beta)=\frac{(\beta^{\tau^{-1}})^3((\beta^{\tau^{-1}})^3+216)^3}{((\beta^{\tau^{-1}})^3-27)^3}=\frac{(\alpha^{\tau^{-1}})^3((\alpha^{\tau^{-1}})^3-24)^3}{(\alpha^{\tau^{-1}})^3-27}=j(E_{\alpha^{\tau^{-1}}}).$$ \smallskip
Therefore $E_\beta \cong E_{\alpha^{\tau^{-1}}}$ and there is an isogeny $\phi_1: E_\alpha \rightarrow E_{\alpha^{\tau^{-1}}}$.  Applying the isomorphism 
$\tau^{-i+1}$ to the coefficients gives an isogeny $\phi_i: E_{\alpha^{\tau^{-(i-1)}}} \rightarrow E_{\alpha^{\tau^{-i}}}$, and therefore an isogeny
$$\iota=\phi_n \circ \phi_{n-1} \circ \cdots \circ \phi_1: E_\alpha \rightarrow E_\alpha.$$ \smallskip
This isogeny has degree $deg(\iota)=3^n$, and we claim that $\Phi_{3^n}(j(E_\alpha),j(E_\alpha))=0$, where $\Phi_m(X,Y)=0$ is the modular equation.  We use several facts from [29].  From [29, Prop. 3.5] the $X'$-coordinate of $\phi(X,Y)$ on the curve $E_\beta$ is given by
$$X'=\frac{-\beta^2}{9\alpha^2} \frac{3X^3+\alpha^2 X^2 +3\alpha X +3}{X^2},$$ \smallskip
and so the $3$-torsion points $\displaystyle P=\left( \frac{-3\beta}{\alpha(\beta-3)}, \frac{\beta-3\omega^i}{\beta-3} \right)$ on $E_\alpha$ map to $(0,0), (0, -1)$ on $E_\beta$.  Further, using [29, Prop. 3.10] the isomorphism between $E_\beta$ and $E_{\alpha^{\tau^{-1}}}$ may be chosen so that the $X_1$-coordinate on $E_{\alpha^{\tau^{-1}}}$ is given in terms of $X'$ on $E_\beta$ by
$$X_1=-\frac{\gamma_1}{\gamma_0} X' + \gamma_1, \quad \gamma_0=\frac{-3\alpha}{\beta(\alpha-3)}, \quad \gamma_1=\left(\frac{-3\beta}{\alpha(\beta-3)} \right)^{\tau^{-1}}. $$
Consequently we have $\phi_1(P)=\pm P^{\tau^{-1}}$, assuming that $\tau$ is defined so that it fixes $\omega$.  Now, successive isogenies $\phi_i$ are defined by conjugating the coefficients in $\phi_1$ by powers of $\tau^{-1}$, so plugging in conjugate points gives
$$\phi_i(P^{\tau^{-(i-1)}})=\pm P^{\tau^{-i}}, \ 1 \le i \le n, \hspace{.07 in} \Longrightarrow \hspace{.07 in} \phi_n(P^{\tau^{-(n-1)}})=\pm P^{\tau^{-n}}=\pm P.$$ \smallskip
Hence, we have $\iota(P)=\pm P$, which implies that $P$ does not lie in the kernel of $\iota$.  From this and the fact that $P \in E_\alpha[3]$ we conclude that $ker(\iota)$ is cyclic, and therefore
$$\Phi_{3^n}(j(E_\alpha),j(E_\alpha))=0,$$
as claimed.  Now by a classical result [10, p.287] we have the factorization
$$\Phi_{3^n}(x,x)=c_n \prod_{-d}{H_{-d}(x)^{r(d,3^n)}},$$
where the product is over discriminants of orders $\textsf{R}_{-d}$ of imaginary quadratic fields and
$$r(d,m)=|\{\xi \in \textsf{R}_{-d}: \xi \ \textrm{primitive}, \ N(\xi)=m \}/\textsf{R}_{-d}^\times|.$$ \smallskip
The exponent $r(d,3^n)$ can only be nonzero when $4^k\cdot 3^n=x^2+dy^2$ has a primitive solution.  The fact that $\textsf{P}_n(x)$ splits in $\textsf{K}_3$ implies that all the conjugate fields of $\mathbb{Q}(\alpha)=\mathbb{Q}(\zeta) \subset \mathbb{Q}_3(\zeta)$ over $\mathbb{Q}$ are unramified at $p=3$.  This implies that $j(E_\alpha)$ is a root of $H_{-d}(x)$ for some $d$ which is not divisible by $3$; hence, $(3,xyd)=1$ and $-d \equiv 1$ (mod 3).  \medskip

Since $j(E_\beta)=j(E_\alpha)^{\tau^{-1}}$ by the first part of the proof, we know that $j(E_\beta)$ is also a root of $H_{-d}(x)$ and by (4.24), the minimal polynomial $p(x)$ of $\alpha$ over $\mathbb{Q}$ is a factor of the polynomial $G_d(x)$ in (3.1).  By the arguments in the proof of Proposition 3.1, $G_d(x)$ factors over $K_1=\mathbb{Q}(\omega)$ as a product of 6 factors, one of which is $p_d(x)$, whose stabilizer in the group $G_{12}$ is generated by the involution $\sigma_1(x)=\frac{3(x+6)}{x-3}$.  Furthermore, we have that
$$(x-3)^{2h(-d)}p_d( \sigma_1(x))=3^{3h(-d)}p_d(x),\eqno{(4.25)}$$ \smallskip
\noindent by Proposition 3.5.  Two factors of $G_d(x)$ are the monic factors
$$\omega^{h(-d)} p_d(\omega x) \ \textrm{and} \ \omega^{2h(-d)} p_d(\omega^2 x),$$
which lie in $K_1[x]$ but not in $\mathbb{Q}[x]$.  (In case $3|h(-d)$, $p_d(x)$ is not a polynomial in $x^3$ because the cube of a root of $p_d(x)$ generates the corresponding ring class field $\Omega_f$ over $\mathbb{Q}$ and therefore cannot lie in a proper subfield of $\Omega_f$.)  It follows that the roots of these two polynomials generate extensions which are ramified at the prime 3. \medskip

Next, the factor $\tilde p(x)=p_d(3\omega)^{-1} p_d^{\sigma_3}(x)$ with $\sigma_3(x)=\frac{3\omega(x+6\omega)}{x-3\omega}$ lies in $\mathbb{Q}[x]$ (see (3.2)).  This is because the mapping $\sigma_4(x)=\frac{3\omega^2(x+6\omega^2)}{x-3\omega^2}=\sigma_1(\sigma_3(x))$, so that (4.25) implies:

$$p_d(3\omega^2)^{-1}  p_d^{\sigma_4}(x)=p_d(3\omega^2)^{-1}  (x-3\omega^2)^{2h(-d)}p_d(\sigma_1 \circ \sigma_3(x))$$
$$=p_d(3\omega^2)^{-1}  (x-3\omega^2)^{2h(-d)} (\sigma_3(x)-3)^{-2h(-d)} 3^{3h(-d)} p_d(\sigma_3(x))$$
$$=p_d(3\omega^2)^{-1}  \frac{(x-3\omega)^{2h(-d)}}{(3(\omega-1))^{2h(-d)}}  3^{3h(-d)} p_d(\sigma_3(x))$$
$$=p_d(3\omega^2)^{-1} (-\omega)^{-h(-d)} (x-3\omega)^{2h(-d)}p_d(\sigma_3(x)).$$ \smallskip

\noindent Putting $x=3\omega$ in (4.25) and noting $\sigma_1(3\omega)=3\omega^2$ gives $(-\omega)^{h(-d)} p_d(3\omega^2)=p_d(3\omega)$, yielding
$$p_d(3\omega^2)^{-1}  p_d^{\sigma_4}(x)=p_d(3\omega)^{-1}  p_d^{\sigma_3}(x),$$ \smallskip
as claimed.  Now, because $\sigma_3$ has order $2$, the roots of $\tilde p(x)$ are the numbers $\sigma_3(\xi)=\frac{3\omega(\xi+6\omega)}{\xi-3\omega}$, where $\xi=3+\pi^3$ runs through the roots of $p_d(x)$ and $\pi$ is a generator in $\Omega_f$ of the ideal $\wp_3$ or its conjugate $\wp_3'$.  Let $\mathfrak{p}$ be a prime divisor of 3 in the field $\Omega_f(\omega)$. If $\mathfrak{p} \nmid {\pi}$, then certainly $\mathfrak{p}$ does not divide the numerator or denominator of
$$\frac{\pi^3+3(1+2\omega)}{\pi^3+3(1-\omega)}=\frac{\sigma_3(\xi)}{3\omega}.$$
On the other hand, if $\mathfrak{p} | \pi$ then $\mathfrak{p}^2 || \pi$, while $\mathfrak{p}^3 || 3(1-\omega)$, and $3(1+2\omega)=3\omega(1-\omega)$.  Hence, $w_\mathfrak{p}(\frac{\pi^3+3(1+2\omega)}{\pi^3+3(1-\omega)})=0$ in any case, and consequently $\sigma_3(\xi) \equiv 0$ (mod 3).  This implies that
$$\tilde p(x) \equiv x^{2h(-d)} \quad (\textrm{mod} \hspace{.05 in} 3)$$ \smallskip
\noindent and proves that $\alpha$ cannot be a root of $\tilde p(x)$.  The same argument obviously applies to the factors $\tilde p(\omega x)$ and $\tilde p(\omega^2 x)$, which are the last two irreducible factors of $G_d(x)$.  This proves that $\alpha$ can only be a root of the factor $p_d(x)$.  Then $\beta=\sigma_1(\alpha^\tau)=\sigma_1(T(\zeta))$ is a root of the same factor.  This shows that every root of $\textsf{P}_n(x)$ is a root of some $p_d(x)$.   \medskip

To complete the proof, we must show that the Frobenius automorphism
$$\tau_d=(\Omega_f/\mathbb{Q}(\sqrt{-d}),\wp_3)$$
has order $n$.  Since $p_d(x)$ divides $\textsf{P}_n(x)$, $\alpha$ is a periodic point of the map $T$ of exact period $n$. The assertion is now clear from (4.9).  \medskip

Finally, since $\textsf{P}_n(x)$ has no multiple roots, and $p_d(x)$ has degree $2h(-d)$, the formula (4.23) is a consequence of what we have proved, Proposition 3.1a), and the fact that
$$\sum_{-d \equiv 1 (3), \ ord(\tau)=n}{2h(-d)}=deg(\textsf{P}_n(x))=2 \sum_{k|n}{\mu (n/k) \ 3^k}, \quad n >1.$$
Conversely, (4.23) implies that every root of $P_n(x)$ is a root of some $p_d(x)$ for which $\tau_d$ has order $n$.  $\square$  \bigskip

This theorem has several interesting consequences.  \bigskip

\noindent {\bf Corollary 1.}  {\it Every periodic point of the maps $T$ or $S$ in the respective sets $\{z: |z|_3 \ge 1\}$ or $\textsf{D}_3-\{3 \}$ is a root of $p_d(x)$ for some integer $d \equiv 2$ (mod $3$), and generates a ring class field over the imaginary quadratic field $K=\mathbb{Q}(\sqrt{-d})$.} \bigskip

To derive the second consequence, we note the following.  From (2.5) we have 
$$j_\alpha=\frac{\alpha^3(\alpha^3-24)^3}{\alpha^3-27} \equiv \alpha^9 \hspace{.05 in} (\textrm{mod} \ \wp_3),$$
if $(\alpha, \beta)$ is the solution of $Fer_3$ that we constructed in Section 2. It follows that $j_\alpha$ is conjugate to $\alpha$ (mod $\mathfrak{p}$) for every prime divisor $\mathfrak{p}$ of $\wp_3$ in $\Omega_f$.  Since $3$ does not divide the discriminant of the class equation $H_{-d}(x)$ and $\Omega_f=K(j_\alpha)$, $H_{-d}(x)$ factors (mod $\wp_3$) into a product of $r=r(d)=h(-d)/n$ distinct polynomials $g_i(x)$ of degree $n$, where $n=\textrm{ord}(\tau)$ is the degree of the prime divisors $\mathfrak{p}_i=(3,g_i(j_\alpha))$ over $\wp_3$ (also the order of $\wp_3$ in the ring class group (mod $f$) of $K$).  Since $\alpha$ is conjugate to $j_\alpha$ (mod $\mathfrak{p}_i$) for each $\mathfrak{p}_i$, its minimal polynomial $m_d(x)$ over $K$ factors exactly the same way, and hence
$$m_d(x) \equiv \prod_{i=1}^r{g_i(x)} \equiv H_{-d}(x) \hspace{.05 in}  (\textrm{mod} \ \wp_3).\eqno{(4.26)}$$ \smallskip
Furthermore, $\beta^\sigma \equiv 0$ (mod $\wp_3$), for all $\sigma \in \textrm{Gal}(\Omega_f/K)$, which implies the congruence
$$p_d(x) \equiv x^{h(-d)} H_{-d}(x) \ (\textrm{mod} \ 3).\eqno{(4.27)}$$

Now let $\mathfrak{D}_n$ denote the set of discriminants $-d \equiv 1$ (mod 3) of orders in imaginary quadratic fields $K=\mathbb{Q}(\sqrt{-d})$ for which the Frobenius automorphism $\tau$ in the corresponding ring class field $\Omega_f$ has order $n$.  Then (4.26) yields a map from the set $\mathfrak{D}_n$ to the power set of the set of all monic irreducible polynomials ($\neq x$) of degree $n$ in $\mathbb{F}_3[x]$:
$$-d \in \mathfrak{D}_n \rightarrow S_d=\{\textrm{irred.} \ g_i(x) \in \mathbb{F}_3[x]: g_i(x) \ | \ H_{-d}(x) \hspace{.05 in} (\textrm{mod} \ 3), \ \textrm{deg}(g_i(x))=n \}.$$ \smallskip
\noindent Different integers $-d \in \mathfrak{D}_n$ yield disjoint sets $S_d$, for the following reason.  If some polynomial $g(x) \in \mathbb{F}_3[x]$ different from $x$ divides both $H_{-d_1}(x)$ and $H_{-d_2}(x)$ (mod 3), then $g(x)$ would divide both $p_{d_1}(x)$ and $p_{d_2}(x)$ (mod $3$), by (4.27).  But then $g(x)^2$ divides $R_n(x)$ (mod $3$) by Theorem 4.2 and (4.19a), which contradicts (4.18).  It follows that
$$\sum_{-d \in \mathfrak{D}_n}{r(d)} \le N_3(n) = \frac{1}{n} \sum_{k|n}{\mu (n/k)3^k}$$ \smallskip
and the above theorem implies that this inequality is actually an equality.  This gives the following corollary. \bigskip

\noindent {\bf Corollary 2.} {\it Every monic irreducible polynomial $f(x) \neq x$ of degree $n \ge 1$ in $\mathbb{F}_3[x]$ divides a unique class polynomial $H_{-d}(x)$ (mod 3) with $-d \in \mathfrak{D}_n$.}  \bigskip

For example, when $n=3$, each of the eight irreducible cubics over $\mathbb{F}_3$ corresponds to a unique discriminant $-d \in \mathfrak{D}_3$, as follows:
$$S_{23}=\{x^3+2x^2+2x+2\}, \quad S_{59}=\{x^3+x^2+x+2\}, \quad S_{83}=\{x^3+2x+2\},$$
$$S_{107}=\{x^3+x^2+2\}, \quad S_{44}=\{x^3+2x^2+1\}, \quad S_{92}=\{x^3+2x^2+x+1\},$$
$$S_{104}=\{x^3+2x+1, x^3+x^2+2x+1\}.$$

This leads to the following question.  The polynomials $g_i(x) \in \mathbb{F}_3[x]$ which divide a given polynomial $p_d(x)$ (mod $3$) in (4.26) appear to satisfy $g_i(x) \equiv g_j(x)$ (mod $x^2$) in $\mathbb{F}_3[x]$; i.e., the $g_i(x)$ all have the same linear and constant terms $ax+b$.  If this does hold for the irreducible factors $g_i(x) \not \equiv x$ of $p_d(x)$ (mod $3$), then there can be no more than $3^{n-2}$ such factors, from which it would follow that $h(-d) \le n \cdot 3^{n-2}$.  If $g(x)$ is a given irreducible of degree $n$ in $\mathbb{F}_3[x]$, is it possible to determine its companions $g_i(x)$ in (4.26) directly from $g(x)$ itself, and thereby determine the class number $h(-d)$ corresponding to the unique polynomial $p_d(x)$ which it divides (mod $3$)?  \medskip

\noindent {\bf Remark.}  The relation (4.23) is equivalent to a class number relation discovered by Deuring [12], [13], which for the prime $p=3$ can be stated as follows:
$$\sum_{-d_{3^n}}{h(-d_{3^n})} =3^n-1,$$
\noindent the sum being taken over all discriminants $-d_{3^n}$ of binary quadratic forms for which the principal form of discriminant $-d_{3^n}$ properly represents $3^n$.  This formula follows by summing (4.23) over the divisors $d \neq 1$ of $n$, and noting from (4.13) that $h(-8)+h(-11)=3-1$, which is the formula corresponding to (4.23) for $n=1$.  The proof given above shows how this formula follows from the theory of the $3$-adic periodic points of the functions $T(z)$ and $S(z)$, which is closely related to the cubic Fermat equation.  Note also that Corollary 2 follows from Deuring's lifting theorem [12].  \bigskip

The above analysis also implies the following theorem.  \bigskip

\noindent {\bf Theorem 4.3.} {\it The unique unramified extension of $\mathbb{Q}_3$ of degree $n \ge 1$ is the splitting field of any polynomial $p_d(x)$ over $\mathbb{Q}_3$ for which $-d \in \mathfrak{D}_n$.  In particular, the maximal unramified extension $\textsf{K}_3$ of $\mathbb{Q}_3$ inside its algebraic closure is the field generated over $\mathbb{Q}_3$ by the roots of all of the polynomials $p_d(x)$, as $-d$ ranges over discriminants $\equiv 1$ (mod $3$) of orders in imaginary quadratic fields.  In other words, $\textsf{K}_3$ is generated over $\mathbb{Q}_3$ by solutions of the cubic Fermat equation $Fer_3$.}  \bigskip

One only has to note that by Corollary 2, there is, for every $n \ge 1$, a corresponding polynomial $p_d(x)$.  Thus there exists at least one discriminant $-d \in \mathfrak{D}_n$. \medskip

Finally, note that Theorem 4.2, together with (4.19a), actually shows that the only numbers $a$ in either of the fields $\bar \mathbb{Q}_3$ or $\mathbb{C}$ which satisfy (4.17) are $a=3$ and the roots of the polynomials $p_d(x)$.  Thus if we define a periodic point of $T(z)$ in $\bar \mathbb{Q}_3$ or $\mathbb{C}$ to be any number $a$ for which there exist numbers $a_1, a_2, \cdots, a_{n-1}$ in $\bar \mathbb{Q}_3$ resp. $\mathbb{C}$ satisfying (4.17), then we have shown the following. \bigskip

\noindent {\bf Theorem 4.4.} {\it The set of periodic points (as defined above) of the multi-valued function $T(z)$ on either of the fields $\mathcal{K}=\bar \mathbb{Q}_3$ or $\mathcal{K}=\mathbb{C}$ coincides with the set}
$$\mathcal{S(K)}=\{3\} \cup \{\alpha \in \mathcal{K}: (\exists n \ge 1)(\exists \ -d \in \mathfrak{D}_n ) \ \textrm{s.t.} \ p_d(\alpha)=0 \}.$$

The same result holds for the algebraic closure $\mathcal{K}=\bar \mathbb{F}_p$ of $\mathbb{F}_p$.

\section{Pre-periodic points of $T(z)$.}

As above, we consider the algebraic function
$$T(z)=\frac{z^2}{3}(z^3-27)^{1/3}+\frac{z}{3}(z^3-27)^{2/3}+\frac{z^3}{3}-6,$$
initially defined on the subset $\{z: |z|_3 \ge 1\} \subset \textsf{K}_3$, where $\textsf{K}_3$ is the maximal unramified, algebraic extension of the $3$-adic field $\mathbb{Q}_3$.  In this section, we consider this function on the algebraic closure $\overline{\mathbb{Q}} \subset \mathbb{C}$ of the rational field $\mathbb{Q}$, but the results apply equally to the algebraic closure $\overline{\mathbb{Q}}_3$ of $\mathbb{Q}_3$.  Denote either of these fields by $\mathcal{K}$.  \medskip

We consider $T(z)$ to be a multi-valued function, and define a pre-periodic point $\xi \in \mathcal{K}$ to be a number for which
$$T^k(\xi)=\alpha, \ \ k \ge 1, \ k \ \textrm{minimal},$$
where $\alpha$ is a periodic point of $T(z)$, i.e., either $3$ or a root of one of the polynomials $p_d(x)$.  The only pre-periodic points for which $T^k(\xi)=3$ are $\xi=3\omega, 3\omega^2$, and then $k=1$.  We shall leave these points out of consideration for the rest of our discussion.  \medskip

Using the same definition that was given in Section 4 above, and setting
$$g(x,y)=(y^2+3y+9)x^3-(y+6)^3,$$
there is a sequence $\xi=\xi_k, \xi_{k-1}, \cdots, \xi_1$ of elements of $\mathcal{K}$ for which
$$g(\xi,\xi_{k-1})=g(\xi_{k-1},\xi_{k-2})= \cdots = g(\xi_2,\xi_1)=g(\xi_1,\alpha)=0.$$
This holds because any branch of the function $T(z)$ satisfies $g(z,T(z))=0$.  We say $\xi_j$ is a pre-periodic point of {\it level} $j$.  \medskip

Since $\alpha$ is a periodic point with some minimal period $n$, then $\alpha_1=T^{n-1}(\alpha)$ satisfies $T(\alpha_1)=\alpha$.  Thus $g(\alpha_1,\alpha)=0$.  However, the form of $g(x,y)$ shows that
$$g(\omega \alpha_1,\alpha)=0, \ \ g(\omega^2 \alpha_1,\alpha)=0,$$
so that $\xi=\omega \alpha_1, \omega^2 \alpha_1$ are pre-periodic with pre-period $k=1$.  It follows, since $\alpha_1$ is also a root of the polynomial $p_d(x)$, that there are at least $4h(-d)$ pre-periodic points having $k=1$ and $p_d(T(\xi))=0$.  On the other hand, the points $\xi$ for which $k=1$ and $T(\xi)$ is a root of $p_d(x)$ are all roots of the resultant
$$R_1(x)=Res_y(p_d(y),g(x,y))=\prod_{p_d(\alpha)=0}{g(x,\alpha)}.$$
This polynomial has degree $6h(-d)$ in $x$, since none of the roots of $p_d(y)$ are roots of $y^2+3y+9$, the coefficient of $x^3$ in $g(x,y)$.  However, $p_d(x) \vert R_1(x)$, so there are at most $6h(-d)-2h(-d)=4h(-d)$ pre-periodic points at level $k=1$.
It follows that there are exactly $4h(-d)$ pre-periodic points at level $k=1$ corresponding to roots of $p_d(x)$, and these are just the roots of the irreducible polynomial (over $\mathbb{Q}$):
$$r_d(x)=p_d(\omega x) p_d(\omega^2 x).$$
Note that any of these pre-periodic points generates the field $\Omega_f(\omega)=\Omega_{3f}$ over $K$. \medskip

Now consider points $\xi$ for which $T^2(\xi)=\alpha$ is a root but $T(\xi)$ is not a root of $p_d(x)$.  Then $T(\xi)=\xi_1$, where $\xi_1=\omega \beta$ or $\omega^2 \beta$ for some root $\beta$ of $p_d(x)$, and $\xi$ is a root of the polynomial
$$R_2(x)=Res_y(r_d(y),g(x,y))=\prod_{p_d(\alpha)=0}{g(x,\omega \alpha)g(x,\omega^2 \alpha)}.\eqno{(5.1)}$$
Recalling that $\sigma_1(z)=\frac{3(z+6)}{z-3}$, we have the following lemma, which is easily verified.  \bigskip

\noindent {\bf Lemma 5.1.}  {\it If $\displaystyle \omega=\frac{-1+\sqrt{-3}}{2}$, then the identity holds:}
$$81 \sqrt{-3} g(x,\omega^2 y) = (y-3)^3 g(x, \omega \sigma_1(y)).$$

From Proposition 3.1 we know that $\sigma_1(z)$ is an involution on the roots of $p_d(x)$.  It follows from (5.1) and this lemma that $R_2(x)$ has at most $6h(-d)$ roots, i.e. $R_2(x)=cs_d(x)^2$ for some constant $c \in \mathbb{Z}$.  We will see that $s_d(x)$ is an irreducible polynomial over $\mathbb{Q}$.  \medskip

Let $w=(k+\sqrt{-d})/2$ or $k+\frac{\sqrt{-d}}{2}$ with $k^2 \equiv d$ (mod $9$) and $k \equiv 1$ (mod 6), as in Theorem 1.  By replacing $k$ by $k+18$ we may assume that $9 || N(w)$, where $N(w)$ is the norm of $w$ to $\mathbb{Q}$.  Then $\{w,9\}$ is a basis for the ideal $\wp_3^2 \cap \textsf{R}_{-d}$ and $j(w/9)$ is a root of the class equation $H_{-d}(x)$ for the discriminant $-d$.  \bigskip

\noindent {\bf Lemma 5.2.} {\it With $\mathfrak{f}(z)$ defined by (2.3), we have $T(\mathfrak{f}(z/3))=\mathfrak{f}(z)$, for $z \in \mathbb{H}$; and more generally, $T^r(\mathfrak{f}(z/3^r))=\mathfrak{f}(z)$, for $r \ge 1$.} \medskip

\noindent {\it Proof.} Using the transformation (4.0) and the fact that $(\mathfrak{f}(z),\mathfrak{g}(z))$ satisfies the equation for $Fer_3$, we see that the pair of functions
$$\left(\frac{9\mathfrak{g}(z)}{\mathfrak{f}(z)(\mathfrak{g}(z)-3)}, \frac{9\mathfrak{g}(z)}{\mathfrak{g}(z)-3}\right)\eqno{(5.2)}$$
satisfies the equation for the curve $E$ in (4.1).  Since
$$\frac{9\mathfrak{g}(z)}{\mathfrak{g}(z)-3} = \frac{3(\mathfrak{g}(z)+6)}{\mathfrak{g}(z)-3}+6=\mathfrak{f}(3z)+6,$$
by (2.8), the point in (5.2) can be written as
$$\left(\frac{\mathfrak{f}(3z)+6}{\mathfrak{f}(z)},\mathfrak{f}(3z)+6\right).$$
Now it follows just as in the equation preceding (4.3) that $g(\mathfrak{f}(z),\mathfrak{f}(3z))=0$.  Putting $z/3$ for $z$ gives the assertion.  $\square$ \bigskip

\noindent {\bf Lemma 5.3.} {\it The number $j(w/3^{r+2})$ generates $\Omega_{3^rf}$ over $K$.} \medskip

\noindent {\it Proof.}  If $d$ is odd, the minimal polynomial of $\frac{w}{3^{r+2}}$ is
$$h(x)=3^{2r+2} x^2-3^rk x+\frac{N(w)}{9}, \ \ (d \ \textrm{odd}).$$
The gcd of the coefficients is $1$, since $9 || N(w)$, and the discriminant of $h(x)$ is $-3^{2r}d=3^{2r}f^2 d_K$.  If $d$ is even, replace $k$ in the polynomial $h(x)$ by $2k$ and we get the same conclusion.  This implies the assertion.  See [10, p. 137].  $\square$  \bigskip

\noindent {\bf Theorem 5.4.} {\it With $w$ defined above, the number $\xi_r=\mathfrak{f}(w/3^{r})$ satisfies $T^{r}(\xi_r)=\alpha$, where $\alpha=\mathfrak{f}(w)$, and $K(\xi_r)=\Omega_{3^{r}f}$ for $r \ge 1$.}  \medskip

\noindent {\it Proof.}  The first assertion is clear from Lemma 5.2.  From the fact that $\mathfrak{g}(w/3)=\sigma_1(\mathfrak{f}(w))=\sigma_1(\alpha)=\beta^{\tau^{-1}}$ (by Proposition 3.2) we see that $\mathfrak{g}(w/3)=\beta^{\tau^{-1}}$.  Hence $\mathfrak{f}(w/3)=\omega^j \alpha^{\tau^{-1}}$, for some $j$, since the point $(\mathfrak{f}(w/3), \mathfrak{g}(w/3))$ lies on $Fer_3$.  Suppose that $j=0$, so $\mathfrak{f}(w/3)=\alpha^{\tau^{-1}} \in \Omega_f$. Then $\sigma_1(\mathfrak{f}(w/3))=\mathfrak{g}(w/9) \in \Omega_f$, so that $j(w/27) \in \Omega_f$ by (2.7b).  But this contradicts Lemma 5.3.  Hence $j \neq 0$, so that $\xi_1=\omega^j \alpha^{\tau^{-1}}$ is a root of $r_d(x)$, and $K(\xi_1)=\Omega_{3f}$. \medskip

Assume inductively that $K(\xi_{r-1})=\Omega_{3^{r-1}f}=K(j(w/3^{r+1}))$, for some $r \ge 2$.  Then $g(\xi_r,\xi_{r-1})=0$ implies that $[K(\xi_r,\xi_{r-1}):K(\xi_{r-1})] \le 3$.  As in the argument in the previous paragraph, $\mathfrak{g}(w/3^{r+1})=\sigma_1(\mathfrak{f}(w/3^r))=\sigma_1(\xi_r)$, so that (2.7b) implies $j(w/3^{r+2}) \in K(\xi_r)$.  Hence, we have that
$$K(\xi_{r-1})=K(j(w/3^{r+1})) \subset K(j(w/3^{r+2})) \subseteq K(\xi_r).$$
Since $[\Omega_{3^r f}:\Omega_{3^{r-1}f}]=3$ for $r \ge 2$, this inclusion shows that $[K(\xi_r):K(\xi_{r-1})]=3$ and $K(\xi)=\Omega_{3^r f}$, as desired.  $\square$  \bigskip

It follows from this theorem that there are at least $[\Omega_{3^r f}:K]=2 \cdot 3^{r-1}h(-d)$ conjugates of $\xi_r=\mathfrak{f}(w/3^r)$ over $K$ and therefore at least this many pre-periodic points at level $r \ge 1$.  Combining this with what we found for level $r=2$, we see that there are exactly $6h(-d)$ pre-periodic points at level $r=2$ attached to roots of $p_d(x)$, and they are all conjugates of $\xi_2$ over the quadratic field $K$.  The polynomial $s_d(x) \in \mathbb{Q}[x]$ has degree $6h(-d)$ and is therefore irreducible over $\mathbb{Q}$.  Note that the $j$-invariants in Lemma 5.3 are algebraic integers, so by (2.7a), the numbers $\xi_r$ are also always algebraic integers.  Hence we can assume $s_d(x) \in \mathbb{Z}[x]$.  \medskip

Now let $s_d^{(r)}(x)$ be the minimal polynomial of the pre-periodic point $\xi_r$ of level $r$ attached to the root $\alpha$ of $p_d(x)$, as before, with $s_d^{(2)}(x)=s_d(x)$.  Then $s_d^{(r+1)}(x)$ is a factor of the resultant $Res_y(s_d^{(r)}(y),g(x,y))$.  On the other hand, the degree of this resultant is $3 \cdot \textrm{deg}(s_d^{(r)}(x))$, so we see inductively that
$$Res_y(s_d^{(r)}(y),g(x,y)) = c_{r+1}s_d^{(r+1)}(x), \ \ r \ge 2.$$ \smallskip
Hence, all pre-periodic points $\xi$ of a given level $r \ge 1$, for which $T^r(\xi)$ is a root of $p_d(x)$, are conjugate over $\mathbb{Q}$, and this gives: \bigskip

\noindent {\bf Theorem 5.5.} {\it All pre-periodic points $\xi \ (\neq 3\omega, 3\omega^2)$ of $T(z)$, such that $T^r(\xi)$ is a root of $p_d(x)$, $r \ge 1$, generate ring class fields over $K=\mathbb{Q}(\sqrt{-d})$ with conductors divisible by $3$.  Every ring class field over $K$ is generated by a periodic point or pre-periodic point of the function $T(z)$.}  \bigskip

This theorem, together with Theorem 4.4, verifies the conjecture put forth in the introduction for the prime $p=3$.  See [9] for other connections between ring class fields and iteration.  \medskip

Remarkably, the polynomial $s_d(x)$ has real roots, even though all the roots of $p_d(x)$ are complex!  Why is this true?  The reason is that for some root $\alpha$ of $p_d(x)$, the quantity $\sigma_1(\alpha)$ is the complex conjugate $\bar \alpha$ of $\alpha$, and for this $\alpha$ the following relation holds:
$$\frac{(\omega \alpha+6)^3}{\omega^2\alpha^2+3\omega \alpha+9}=\frac{(\omega^2 \bar \alpha+6)^3}{\omega \bar \alpha^2+3\omega^2 \bar \alpha+9}.$$
Therefore, at least one of the cube roots $\xi$ of this real quantity is real, as well, and this cube root satisfies the equation $g(\xi,\omega \alpha)=0$, so that $T(\xi)=\omega \alpha$.

\section{Aigner's Conjecture and the rank of $E(\Sigma)$.}

In this section and the next we shall give several applications of the results obtained in sections 3 and 4.  We first restrict our attention to the case when $f=1$, so that $-d=d_K=\textrm{disc}(K/\mathbb{Q})$, and $\Omega_f=\Sigma$ is the Hilbert class field of $K$.  In this section we once again work on the curve
$$E: \ \ Y^2-9Y=X^3 -27,$$
for which $E(\mathbb{Q})=\{O,(3,0),(3,9)\}$, so that $(3,0)$ and $(3,9)$ have order $3$ on $E$.  We refer to these three points as the {\it trivial points} on $E$.  \bigskip

\noindent {\bf Lemma 6.1.} {\it Assume that $(\alpha,\beta)$ is a point on $Fer_3$ with $\alpha, \beta \neq 0, 3, \infty$ and $\alpha \neq \beta$.  Let
$$P= (x,y)=\left( \frac{9\beta}{\alpha(\beta-3)},\frac{9\beta}{\beta-3} \right) \ and \ Q= (x_1,y_1)=\left( \frac{9\alpha}{\beta(\alpha-3)},\frac{9\alpha}{\alpha-3} \right)$$
be the points on $E$ corresponding to $(\alpha, \beta)$ and $(\beta,\alpha)$ on $Fer_3$.  Then $P+Q=(3,9)$.}  \medskip

\noindent {\it Proof.} We have the rational identity
$$-3-x-x_1+\left( \frac{y-y_1}{x-x_1} \right)^2 =\frac{3(2\alpha \beta-3\alpha-3\beta)(\alpha^3 \beta^3-27\alpha^3-27\beta^3)}{(\alpha \beta-3\alpha-3\beta)^2 \alpha \beta (\alpha-3)(\beta-3)}.$$ \smallskip
\noindent Since the resultant of $\alpha \beta-3\alpha-3\beta$ and $\alpha^3 \beta^3-27\alpha^3-27\beta^3$ is $-243\alpha^4(\alpha-3)\neq 0$, the denominator in the last expression is not zero, so the right hand side is zero.  Hence, the $X$-coordinate of the sum $P+Q$ on $E$ is $3$.  The fact that the $Y$-coordinate of $P+Q$ is $Y=9$ follows from the equation
$$9-Y=\left( \frac{y-y_1}{x-x_1} \right) (3-x)+y=\frac{-3\alpha \beta}{\alpha \beta-3\alpha-3\beta} \left( \frac{3(\alpha \beta-3\alpha-3\beta)}{\alpha(\beta-3)}\right)+\frac{9\beta}{\beta-3}=0.$$
\noindent This proves the lemma. $\square$  \bigskip

We apply Lemma 6.1 to the point
$$P_d=\left( \frac{9\beta}{\alpha(\beta-3)},\frac{9\beta}{\beta-3} \right)=\left(\frac{\alpha^\tau+6}{\alpha},\alpha^\tau+6 \right)$$
in $E(\Sigma)$ whose coordinates are given in (4.2).  By Proposition 3.2, the map $\phi$ defined by $\alpha^\phi=\sigma_1(\alpha^\tau)=\beta$ is an automorphism of $\Sigma$ of order 2 and therefore switches $\alpha$ and $\beta$.  Note that the map on the quadratic field $K$ induced by $\phi$ also interchanges the ideals $\wp_3$ and $\wp_3'$.  Lemma 6.1 implies that $P_d+P_d^\phi=(3,9)$ and therefore
$$P_d^\phi=(3,9)-P_d.$$ \smallskip
\noindent Letting $G=\textrm{Gal}(\Sigma/K)$, we have therefore that
$$P_d^{\sigma \phi}=P_d^{\phi \sigma^{-1}}=(3,9)-P_d^{\sigma^{-1}}, \ \  \sigma \in G, \eqno{(6.1)}$$
and hence
$$\left( \sum_{\sigma \in G}{P_d^\sigma} \right)^\phi=\sum_{\sigma \in G}{P_d^{\sigma \phi}} = \sum_{\sigma \in G}{(3,9)} - \sum_{\sigma \in G}{P_d^{\sigma^{-1}}}=[h(K)](3,9)-\sum_{\sigma \in G}{P_d^\sigma} ,$$
where $h(K)$ is the class number of $K$.  Hence, $\displaystyle Q_K=\sum_{\sigma \in G}{P_d^\sigma}$ satisfies
$$Q_K^\phi=[h(K)](3,9) - Q_K.\eqno{(6.2)}$$ \smallskip
(Cf. [6, eq. (8)].)  Now if $\displaystyle Q_K$ lies in $E(\mathbb{Q})$, then $Q_K^\phi=Q_K$, and the last equation gives
$$[2]Q_K=[2]\sum_{\sigma \in G}{P_d^\sigma} =[h(K)] (3,9).$$
By the doubling formula on the curve $E$, if $Q_K \neq O$ and $x=x(Q_K)$, then either $4x^3-27=0$ or
$$\frac{x^4+54x}{4x^3-27}=3 \ \Rightarrow \ (x-3)(x^3-9x^2-27x-27)=0.$$
Since both cubics $4x^3-27$ and $x^3-9x^2-27x-27$ are irreducible over $\mathbb{Q}$ (both have roots which generate ramified extensions at $p=3$), it follows that $x(Q_K)=3$ and
$$Q_K \in E(\mathbb{Q}) \ \Rightarrow \ Q_K=[2h(K)] (3,9)=[h(K)](3,0).\eqno{(6.3)}$$ \smallskip
In other words, if rational, $Q_K$ is determined by the residue class of $h(K)$ (mod 3). \bigskip

Now let $\mathfrak{p}$ be a prime divisor in $\Sigma$ of $\wp_3$, and consider the reduction of $E$ mod $\mathfrak{p} =\tilde E$ taking $P$ to $P$ (mod $\mathfrak{p}$) =$\tilde P$.  We have that
$$\tilde E: \quad y^2 =x^3,$$ \smallskip
\noindent with singular point $S=(0,0)$ mod $\mathfrak{p}$.  The set $\tilde E_{ns}$ of non-singular points on $\tilde E$ forms a group isomorphic to $k_\mathfrak{p}^+=(R_{\Sigma}/\mathfrak{p})^+$, which is a vector space of dimension $f(\mathfrak{p}/\wp_3)$ over $\mathbb{F}_3$.  Since the point $S$ is a cusp and the tangent line to $\tilde E$ at $S$ is the line $y=0$, an isomorphism between $\tilde E_{ns}$ and $k_\mathfrak{p}^+$ is given by the map

$$(x,y) \in \tilde E_{ns} \longrightarrow \frac{x}{y}=\frac{1}{\alpha} \in k_\mathfrak{p}^+.\eqno{(6.4)}$$ \smallskip
(See [38, p. 56].)  Furthermore, the set
$$E_\mathfrak{p}(\Sigma)=\{P \in E(\Sigma): \hspace{.05 in} \tilde P \in \tilde E_{ns} \}$$ \smallskip
\noindent is a subgroup of $E(\Sigma)$.  Now the point $P_d=\left( \frac{\alpha^\tau+6}{\alpha}, \alpha^\tau +6 \right)$ has nonzero coordinates (mod $\mathfrak{p}$), since $\alpha$ and therefore also $\alpha^\tau$ has no prime divisors in common with $\wp_3$, by Theorem 3.4a) and the fact that $\tau=(\Sigma/K,\wp_3)$ fixes $\mathfrak{p}$.  Therefore, $\tilde P_d \neq S$ and hence $P_d \in E_\mathfrak{p}(\Sigma)$.  Since the automorphisms $\sigma \in G$ permute the prime divisors of $\wp_3$ among themselves, we also have $P_d^\sigma \in E_\mathfrak{p}(\Sigma)$, for $\sigma \in G$.  Since $E_\mathfrak{p}(\Sigma)$ is a subgroup, the sum
$$Q_K =\sum_{\sigma \in G} {P_d^\sigma} \in E_\mathfrak{p}(\Sigma).\eqno{(6.5)}$$ \smallskip
\noindent But the reductions of the points $(3,0)$ and $(3,9)$ are $S$, so that $(3,0), (3,9) \notin E_\mathfrak{p}(\Sigma)$.  Therefore $Q_K \neq (3,0)$ or $(3,9)$.  Now (6.3) shows that $Q_K \notin E(\mathbb{Q})$ if the class number of $K$ is prime to 3.  \bigskip

\noindent {\bf Theorem 6.2.} {\it Aigner's conjecture is true for the imaginary quadratic field $K$ if $\textrm{disc}(K/\mathbb{Q})=-d \equiv 1$ (mod 3) and  $3$ does not divide the class number $h(K)$.  In this case, the point $Q_K$ has infinite order in $E(K)$.} \medskip

\noindent {\it Proof.} If $n \ge 2$ is a natural number for which $(3,n)=1$, then from (6.2) we have
$$[n]Q_K^\phi=[nh(K)](3,9) - [n]Q_K,$$
and this gives immediately that $[n]Q_K \neq O$.  On the other hand, if $3^k n$ is the order of $Q_K$ in $E(K)$, for some $ k \ge 1$, then $[3^{k-1}n]Q_K$ is a point of order $3$.  But the only points of order $3$ in $E(K)$ are $(3,0)$ and $(3,9)$, since $E[3]$ is elementary abelian of order $9$ and contains the points $(3\omega,0), (3\omega, 9)$. Thus, $[3^{k-1}n]Q_K\in \{(3,0), (3,9)\}$.  But this is impossible, by the same argument as above, since $[3^{k-1}n]Q_K \in E_\mathfrak{p}(\Sigma)$.  Thus, $Q_K$ cannot have finite order.  $\square$ \bigskip

It is known that the relative density of quadratic fields $K=\mathbb{Q}(\sqrt{-d})$ with $-d \equiv 1$ (mod 3) and $h(-d) \not \equiv 0$ (mod 3) is at least $1/2$.  (See [33, Thm. 1] and [8, Lemma 2.2].)  In fact, there are 229 square-free integers $d_0 \equiv 2$ (mod $3$) less than $1000$, and if the discriminant of $\mathbb{Q}(\sqrt{-d_0})$ is $-d$, the class number $h(-d)$ is not divisible by $3$ for 165 of these $d_0$.  Thus, $3$ divides $h(-d)$ for only $64/229 \approx .279$ of these discriminants.   \medskip

This result is a counterpart to the results of Fueter and Aigner for the first and second families of quadratic fields mentioned in the Introduction, since those results also require the class number of the imaginary quadratic field in question to be relatively prime to $3$.  Equation (6.3) shows furthermore that when $h(K) \equiv 0$ (mod $3$), the point $Q_K \in E(K)$ can only be trivial if $Q_K=O$.  \medskip

We can also apply the above argument to intermediate fields between $K$ and $\Sigma$.  If $L$ is any intermediate field between $K$ and $\Sigma$ and $H$ is the corresponding subgroup of $G$ in the Galois correspondence, then the point

$$Q_L=\sum_{\sigma \in H}{P_d^\sigma}$$ \smallskip
\noindent satisfies the identity
$$Q_L^\phi=|H| (3,9)-Q_L.\eqno{(6.6)}$$ \smallskip

\noindent The same arguments used to prove Theorem 6.2 give the following result.  \bigskip

\noindent {\bf Theorem 6.3.} {\it Let $H$ be a subgroup of $G=\textrm{Gal}(\Sigma/K)$ whose order is not divisible by 3, and let $L$ be the subfield of $\Sigma$ corresponding to $H$ in the sense of Galois theory.  Then the cubic Fermat equation has a nontrivial solution in the field $L$. Such a solution corresponds to the point $Q_L$ in $E(L)$ defined by $Q_L=\sum_{\sigma \in H} P_d^\sigma$.}  \bigskip

This theorem allows us to prove the following result. \bigskip

\noindent {\bf Theorem 6.4.} {\it Let $K=\mathbb{Q}(\sqrt{-d})$ with $-d \equiv 1$ (mod 3) and $3 \mid h(K)$, and let $d_1$ be the square-free part of $d$.  Then there are infinitely many positive fundamental discriminants $D \equiv 1$ (mod $12d_1$) for which the cubic Fermat equation has a nontrivial solution in the quartic field $L=\mathbb{Q}(\sqrt{-d}, \sqrt{D})$.} \medskip

For the proof we require a lemma from Byeon's paper [8]. \bigskip

\noindent {\bf Lemma.} ([8, Prop. 3.1] ) {\it Let $d_1$ be the square-free part of $d$.  For any square-free integer $t$ there are infinitely many positive fundamental discriminants $D \equiv 1$ (mod $12d_1$) (with positive density) for which the class numbers of the fields $K_1=\mathbb{Q}(\sqrt{tD})$ and $K_2=\mathbb{Q}(\sqrt{D})$ are not divisible by $3$.} \bigskip

\noindent {\it Proof of Theorem 6.4.} Take $t$ in the lemma to be the square-free part of $d_K=-d$.  Then for infinitely many fundamental discriminants $D \equiv 1$ (mod $12d_1$) we have $h(-dD) \not \equiv 0$ (mod $3$).  Hence, Theorem 6.3 applies to the field $L=\mathbb{Q}(\sqrt{-d}, \sqrt{D})$, which is an intermediate field between $K_1=\mathbb{Q}(\sqrt{-dD})$ and its genus field and is therefore contained in the Hilbert class field of $K_1$.  We conclude that the point $Q_L$ is a nontrivial point in $E(L)$.  $\square$ \bigskip

Note that under the isomorphism between $\tilde E_{ns}$ and $k_\mathfrak{p}^+$ (for $\mathfrak{p}  |  \wp_3$) the point $\tilde Q_L = \sum_{\sigma \in H}{\tilde P_d^\sigma}$ maps to the residue class in $k_\mathfrak{p}^+$ of the element
$$ \sum_{\sigma \in H}{\frac{1}{\alpha^\sigma}} = Tr_{\Sigma/L}\left( \frac{1}{\alpha} \right) (\textrm{mod} \ \mathfrak{p}). \eqno{(6.7)}$$
Thus $\tilde Q_L = \tilde O$ if and only if $\mathfrak{p} \mid Tr_{\Sigma/L}(1/\alpha)$. This gives the following criterion in the case that $h(K)$ is divisible by $3$.  \bigskip

\noindent {\bf Theorem 6.5.} {\it If $3 \mid h(K)$ and the solution $(\alpha,\beta)$ of $Fer_3$ in the Hilbert class field $\Sigma$ of $K$ satisfies the condition
$$Tr_{\Sigma/K}\left(\frac{1}{\alpha}\right) \not \equiv 0 \ (\textrm{mod} \ \wp_3) \ \ \textrm{in} \ K,$$
then $Q_K$ is a nontrivial point in $E(K)$.}  \bigskip

From (4.26) the minimal polynomial $m_d(x)$ of $\alpha$ over $K$ satisfies the congruence
$$m_d(x) \equiv H_{-d}(x) \  (\textrm{mod} \ \wp_3).$$ \smallskip
Now if $m_d(x)=x^h+\cdots+c_1x+c_0$, we have that $Tr_{\Sigma/K}(1/\alpha)=-c_1/c_0$.  Since 3 does not divide $H_{-d}(0)$ (see the proof of Proposition 3.2), the above congruence gives
$$Tr_{\Sigma/K}(1/\alpha)=-\frac{c_1}{c_0} \equiv -\frac{H_{-d}'(0)}{H_{-d}(0)} \ (\textrm{mod} \ \wp_3).$$
Using (4.27), this yields the following restatement of Theorem 6.5.
\bigskip

\noindent {\bf Corollary of Theorem 6.5.} {\it If $3$ divides $h(-d)$ but does not divide $H_{-d}'(0)$, the cubic Fermat equation has a nontrivial solution in $K=\mathbb{Q}(\sqrt{-d})$. The condition $H_{-d}'(0) \not \equiv 0$ (mod $3$) is equivalent to the coefficient of $x^{h(-d)+1}$ in $p_d(x)$ not being divisible by $3$.}  \bigskip

For example, the discirminants $-23, -59, -83$ satisfy the hypothesis of the corollary because $h(-d)=3$ for each of these discriminants and
\begin{eqnarray*}
p_{23}(x) &=& x^6+11x^5+65x^4+191x^3+441x^2+405x+675,\\
p_{59}(x) &=& x^6+22x^5+208x^4-40x^3+144x^2-3456x+6912,\\
p_{83}(x) &=& x^6 +6x^5+560x^4-1384x^3+576x^2-12960x+43200.
\end{eqnarray*}
Thus, there is a nontrivial solution of $Fer_3$ in each of the fields $\mathbb{Q}(\sqrt{-23})$, $\mathbb{Q}(\sqrt{-59})$, and $\mathbb{Q}(\sqrt{-83})$.  In fact, we have
$$Q_K=\cases{(1, \frac{9+\sqrt{-23}}{2}), &if \ $d=23$, \cr (-2, \frac{9+\sqrt{-59}}{2}), &if \ $d = 59$, \cr (-8, \frac{9+5\sqrt{-83}}{2}), &if \ $d=83$.}$$

The discriminant $-d=-107$ does not satisfy this condition, since
$$p_{107}(x)=x^6-74x^5+1680x^4-6184x^3+2736x^2-43200x+172800;$$
nonetheless, there is still a nontrivial solution of $Fer_3$ in $K=\mathbb{Q}(\sqrt{-107})$, namely
$$Q_K=\left(\frac{-26}{9}, \frac{9}{2}+\frac{29}{54} \sqrt{-107}\right).$$
We will see below that $Q_K$ can be deduced to be nontrivial in this case by the criterion of Theorem 7.1 (also see Proposition 7.2).  These computations and the factorization of the iterated resultant $R_3(x)$ in Section 4 show that $Fer_3$ has a nontrivial solution in all of the fields $K=\mathbb{Q}(\sqrt{-d})$ for which $-d \equiv 1$ (mod $3$) and $h(K)=3$.  \medskip

We also note that $Q_K$ is sometimes the zero point in $E(K)$.  It is known that the quadratic twist of $X_0(27)\cong Fer_3$ by $K=\mathbb{Q}(\sqrt{-5219})=\mathbb{Q}(\sqrt{-17 \cdot 307})$ has rank $3$.  Assuming the Birch-Swinnerton-Dyer conjecture for this curve, the Gross-Zagier theorem [20] on the derivatives of $L$-functions of elliptic curves would imply that the canonical height of the trace to $K$ of a Heegner point on $Fer_3(\Sigma)$ is necessarily zero.  This would imply that $Q_K$ has finite order for this field $K$, and therefore $Q_K=O$.  (I am grateful to the referee for these remarks.)  The coefficient of $x^{25}=x^{h(-5219)+1}$ in $p_{5219}(x)$ turns out to be
$$-10952 776646 76630 58811 93496 594488 18230 202411 02307328$$
$$= -2^{26}3^2(107071)(29757069131)(56916714418935524735887103),$$ 
so that the conditions of Theorem 6.5 and its corollary are not satisfied; neither is (1.4) of the Introduction, since that condition holds if and only if $9$ does not divide the coefficient of $x^{h(-d)+1}$ in $p_{d}(x)$. See Proposition 7.2 below.  \medskip

We finish this section by proving the following theorem, in which we consider the solutions $(\alpha,\beta)$ of $Fer_3$ in ring class fields $\Omega_f$ with $(f,3)=1$.  \bigskip

\noindent {\bf Theorem 6.6.} {\it Let $\mathfrak{p}$ be a prime divisor of $\wp_3$ in the ring class field $\Omega_f$ of $K$ (with $(f,3)=1$) and let
$$\ell=\textrm{dim}_{\mathbb{F}_3} \langle \frac{1}{\alpha^\sigma} (\textrm{mod} \ \mathfrak{p}), \ \sigma \in G\rangle, \ G= \textrm{Gal}(\Omega_f/K)$$
be the dimension of the vector space generated by the residue classes of the numbers $1/\alpha^\sigma$ in  $R_{\Omega_f}/\mathfrak{p}=\mathbb{F}_{3^n}$, where $\alpha$ is the number defined in Theorem 1.  Then the rank of $Fer_3$ over the ring class field $\Omega_f$ is at least $\ell$. In particular, the rank of $Fer_3$ over the Hilbert class field $\Sigma$ is always at least $1$.}  \smallskip

\noindent {\it Proof.} Let $\{1/\alpha^\sigma | \ \sigma \in S \subseteq G\}$ be a maximal set of linearly independent conjugates of $1/\alpha$ over $\mathbb{F}_3$.  we claim that the corresponding set $\{P_d^\sigma | \ \sigma \in S\}$ is a set of independent points in $E(\Omega_f)$.  Suppose that $\sum_{\sigma \in S}{c_\sigma P_d^\sigma}=O$, for some integers $c_\sigma$.  We then have
$$\sum_{\sigma \in S}{c_\sigma \tilde P_d^\sigma}=\tilde O \ \ \textrm{in} \ R_{\Omega_f}/\mathfrak{p}.$$
By the isomorphism (6.4), this implies that
$$\sum_{\sigma \in S}{c_\sigma \frac{1}{\alpha^\sigma}} \equiv 0 \ (\textrm{mod} \ \mathfrak{p}) \ \textrm{in} \ R_{\Omega_f}/\mathfrak{p}.$$
It follows that $3 | c_\sigma$ for all $\sigma \in S$.  Hence, the sum
$$\sum_{\sigma \in S}{\left(\frac{c_\sigma}{3}\right) P_d^\sigma} \eqno{(6.8)}$$
is a point of order $1$ or $3$ in $E(\Omega_f)$.  But the only points of order $3$ in $E(\Omega_f)$ are $(3,0)$ and $(3,9)$, by the same argument as in the proof of Theorem 6.2, since $\omega \notin \Omega_f$.  Since the sum in (6.8) lies in $E_\mathfrak{p}(\Omega_f)$, with the same notation as above, it follows that
$$\sum_{\sigma \in S}{\left(\frac{c_\sigma}{3}\right) P_d^\sigma} =O.$$
This argument can be repeated indefinitely if some $c_\sigma \neq 0$, giving a contradiction.  Hence the points $P_d^\sigma$, for $\sigma \in S$, are independent and they generate a subgroup of $E(\Omega_f)$ of rank $\ell$.  $\square$ \bigskip

\noindent {\bf Corollary.}  {\it If $\mathfrak{D}_n$ is defined as in Section 4, then for at least one discriminant $
-d=d_Kf^2 \in \mathfrak{D}_n$, the rank of $Fer_3(\Omega_f)$ is at least $n$. Thus, the rank of $Fer_3(\Omega_f)$ is unbounded over ring class fields $\Omega_f$ of imaginary quadratic fields $K$ of the fourth family.}  \medskip

\noindent {\it Proof.} By Corollary 2 of Theorem 4.2, there is at least one value of $-d \in \mathfrak{D}_n$ and an irreducible $f(x) \in \mathbb{F}_3[x]$ of degree $n$ dividing $H_{-d}(x)$ (mod $3$) for which the reciprocal $1/\alpha$ of a root of $f(x)$ generates a normal basis of $\mathbb{F}_{3^n}$ over $\mathbb{F}_3$.  For this $d$, $\ell=n$ in Theorem 6.5.  This proves the corollary.  $\square$ \bigskip

An easy computation shows that the only {\it cubics} in $\mathbb{F}_3[x]$ for which $1/\alpha$ does {\it not} generate a normal basis over $\mathbb{F}_3$ are $f(x)=x^3+x^2+2, x^3+2x^2+1$.  Hence, the rank of $Fer_3(\Sigma)$ is at least $3$ for each of the discriminants in the set $\{-23,-59, -83, -104\} \subset \mathfrak{D}_3$.  \medskip

For $d=5219$, the polynomials $p_{5219}(x)$ and $q_{5219}(x)$ satisfy
$$p_{5219}(x) \equiv q_{5219}(x) \equiv x^{24}(x^8+2x^7+x^3+1)(x^8+x^5+2x^4+2x^2+1)$$
$$ \ \ \ \ \times (x^8+x^7+2x^6+x^5+x^4+2x^3+x^2+1) \ (\textrm{mod} \ 3).$$
Using the third $8$-th degree polynomial in this factorization, it was computed that the elements $1/\alpha^{\tau^i}$, for $0 \le i \le 6$, are independent over $\mathbb{F}_3$, so $\ell =7$.  Thus, $\textrm{rank}(E(\Sigma)) \ge 7$ for $d=5219$.

\section{The formal group of $E$.}

In this section we will use the series for $T(z)$ in (4.6) to compute a formal group for the elliptic curve $E$.  We replace $z$ by $1/z$ in $T(z)$ and define the following functions:

$$x(z)=zT\left( \frac{1}{z} \right)+6z=\frac{1}{z^2}-\sum_{k=1}^\infty{\frac{3^{2k-1}b_k}{k!} z^{3k-2}},$$
$$y(z)=T\left( \frac{1}{z} \right)+6=\frac{1}{z^3}-\sum_{k=1}^\infty{\frac{3^{2k-1}b_k}{k!} z^{3k-3}}.$$ \smallskip

\noindent For all $z \in \textsf{K}_3$ for which $|z|_3 \le 1$ these series are convergent in $\textsf{K}_3$ and $(x(z),y(z))$ is a point on the curve $E$.  Following [38, Ch. IV.1] but with a change of sign, we define

$$w(z) = \frac{1}{y(z)} =z^3 \left(1+\sum_{k=1}^\infty{3^k c_k z^{3k}} \right), \quad c_k \in \mathbb{Z},$$ 
$$=z^3+9 z^6+135 z^9+2430 z^{12}+48114 z^{15}+1010394 z^{18}+22084326 z^{21}+\cdots,$$ \smallskip
so that $(z,w(z))=(x(z)/y(z), 1/y(z))$ lies on the curve
$$E': \quad w=f(z,w)=z^3+9w^2-27w^3.\eqno{(7.1)}$$
The coefficients in the series for $w(z)$ have the form $3^k c_k$, with $c_k \in \mathbb{Z}$, since $3^{k-1}b_k/k! \in \mathbb{Z}$, by the discussion following (4.6).  These coefficients are quite remarkable, in that they all seem to factor into very small primes.  The largest prime factor of any $c_k$ for $k \le 32$ is $p=89$.   \medskip

We now compute a formal group $F_E(z_1,z_2)$ using the model $E'$.  Setting
$$\lambda=\frac{w(z_2)-w(z_1)}{z_2-z_1}, \quad \nu=w(z_1)-\lambda z_1,$$ \smallskip
we have $\lambda, \nu \in \mathbb{Z}[[z_1,z_2]]$ and the points $P_1=(z_1,w(z_1)), P_2=(z_2,w(z_2))$ lie on the line $w=\lambda z+\nu$.  If $P_1+P_2+P_3=O$ and $P_3=(z_3,w_3)$, then substituting $w=\lambda z+\nu$ in the equation $f(z,w)-w=0$ and considering the coefficient of $z^2$ yields the relation
$$z_1+z_2+z_3=-9\lambda^2 \frac{1-9\nu}{1-27\lambda^3}.$$ \smallskip
If $|z_i|_3 \le 1$ for $i=1,2$, then $|\lambda|_3, |\nu|_3 \le 1$ and the last relation shows that $|z_3|_3 \le 1$ also.  By Lemma 4.1 with $z=1/z_3$ and $w = 1/w_3 -6$ we conclude that $w_3=w(z_3)$.  \medskip

To determine the formal group we still need the $z$-coordinate $i(z)$ of the inverse of the point $P=(z,w(z))$.  Since the inverse of $(x(z),y(z))$ on $E$ is $(x(z), 9-y(z))$, we have
$$i(z)=\frac{x(z)}{9-y(z)}=-\frac{x(z)}{y(z)} \sum_{k=0}^\infty{\left( \frac{9}{y(z)} \right)^k} = -z \sum_{k=0}^\infty{3^{2k} w(z)^k}.$$ \smallskip
Thus, with
$$z_3=-z_1-z_2-9\lambda^2 \frac{1-9\nu}{1-27\lambda^3}=G(z_1,z_2),$$
we have
$$F_E(z_1,z_2)=i(z_3)=-z_3 \sum_{k=0}^\infty{3^{2k} w(z_3)^k}=\left(z_1+z_2+9\lambda^2 \frac{1-9\nu}{1-27\lambda^3}\right)\sum_{k=0}^\infty{3^{2k} w(z_3)^k}.$$ \smallskip

\noindent In particular,
$$F_E(z_1,z_2) \equiv z_1+z_2 \hspace{.07 in} (\textrm{mod} \hspace{.05 in} 9), \quad |z_1|_3, |z_2|_3 \le 1. \eqno{(7.2)}$$
We also have
\begin{eqnarray*}
F_E(z_1,z_2) &\equiv& (z_1+z_2+9\lambda^2)(1+9w(z_3)) \quad (\textrm{mod} \ 3^4)\\
&\equiv& z_1+z_2 +9(-(z_1+z_2)^4+(z_1^2+z_1z_2+z_2^2)^2) \quad (\textrm{mod} \ 3^4)\\
&\equiv& z_1+z_2 -9z_1 z_2 (2z_1^2+3z_1 z_2+2z_2^2) \quad (\textrm{mod} \ 3^4),
\end{eqnarray*}
and therefore
$$F_E(z_1,z_2) \equiv z_1+z_2+9z_1z_2(z_1^2+z_2^2) \ \ (\textrm{mod} \ 27).\eqno{(7.3)}$$

Note that $F_E(z_1,z_2)$ converges for all integral elements $z_1, z_2 \in \textsf{K}_3$, which is a bit stronger than the usual condition for convergence of the formal group.  We now use (7.2) to prove \bigskip

\noindent {\bf Theorem 7.1.} {\it Assume that $3 \mid h(K)$ but that $\wp_3^2$ does not divide $Tr_{\Sigma/K}(1/\alpha)$.  Then $Q_K$ is a nontrivial point in $E(K)$ and Aigner's conjecture holds for the field $K=\mathbb{Q}(\sqrt{-d})$.} \medskip

\noindent {\it Proof.} Let $\mathfrak{p}$ be any prime divisor of $\wp_3$ in $\Sigma$, and consider the embedding of $\Sigma/K$ in the completion $\Sigma_\mathfrak{p}/\mathbb{Q}_3$, so that $\pi=3$ is a prime element for $\mathfrak{p}$.  Let $z=x/y=1/\alpha$ be the $z$-coordinate of the point $P_d=(x,y)$ on the model (7.1).  Then the $z$-coordinate of $P_d^\sigma$ is $z(P_d^\sigma)=1/\alpha^\sigma$ for $\sigma \in G=\textrm{Gal}(\Sigma/K)$, and we have $|z(P_d^\sigma)|_3=1$ (by Theorem 3.4a)).  Therefore, (7.2) implies that the $z$-coordinate of the point $Q_K$ satisfies
$$z(Q_K) \equiv \sum_{\sigma \in G}{z(P_d^\sigma)} = \sum_{\sigma \in G}{\frac{1}{\alpha^\sigma}}=Tr_{\Sigma/K}\left(\frac{1}{\alpha} \right) \quad (\textrm{mod} \hspace{.05 in} \mathfrak{p}^2).$$ 
Now if $Q_K=O$, then $z(Q_K)=0$.  But then the last congruence implies that 
$$Tr_{\Sigma/K}\left(\frac{1}{\alpha} \right) \equiv 0 \ (\textrm{mod} \ \mathfrak{p}^2), \quad \mathfrak{p} \mid \wp_3.$$ \smallskip
Since this holds for all prime divisors $\mathfrak{p}$ of $\wp_3$ in $\Sigma$, we get that $\wp_3^2 \mid Tr_{\Sigma/K}(1/\alpha)$, contrary to hypothesis.  Hence $Q_K \neq O$, so the discussion in Section 6 implies that $Q_K$ is a nontrivial point in $E(K)$.  $\square$   \medskip

\noindent {\bf Corollary 1.} {\it In the situation of Theorem 7.1, the point $Q_K$ has infinite order on $E$.}\medskip

\noindent {\it Proof.} From (7.2), $z([n]Q_K) \equiv n z(Q_K) \equiv n \cdot Tr(1/\alpha)$ (mod $\wp_3^2$), so that $Q_K$ cannot have order $n$, if $3 \nmid n$.  On the other hand, if $3 \mid n$, then $Q_K$ cannot have order $n$ by the proof of Theorem 6.2.  $\square$ \medskip

We can express this theorem using the class equation $H_{-d}(x)$, as follows.  We have

$$j_\alpha=\frac{\alpha^3(\alpha^3-24)^3}{\alpha^3-27} \equiv (\alpha^3 - 6)^3 \ (\textrm{mod} \ \wp_3^2),$$
so
$$j_\alpha-6 \equiv (\alpha^3 - 6)^3 -6 \equiv T^2(\alpha) \ (\textrm{mod} \ \mathfrak{p}^2), \ \mathfrak{p} \mid \wp_3,$$ \smallskip
\noindent by (4.7), under the embedding of $\Sigma$ in the $\mathfrak{p}$-adic completion $\Sigma_\mathfrak{p} \subset \textsf{K}_3$, where $\mathfrak{p}$ is a prime divisor of $\wp_3$ in $\Sigma$.  But from (4.9), $T^2(\alpha)=\alpha^{\tau^2}$, which implies that
$$j_\alpha-6 \equiv \alpha^{\tau^2} \hspace{.07 in} (\textrm{mod} \hspace{.05 in} \wp_3^2),$$ \smallskip
and therefore $H_{-d}(x+6) \equiv m_d(x)$ (mod $\wp_3^2)$ in $K$, where $m_d(x)$ is the minimal polynomial of $\alpha$ over $K$.  Since $\pm N_{\Sigma/K}(\alpha)$, the constant term of $m_d(x)$, is relatively prime to $\wp_3$, it follows that $\wp_3^2$ divides $Tr_{\Sigma/K}(1/\alpha)$ if and only if $9 \mid H_{-d}'(6)$.  \bigskip

\noindent {\bf Corollary 2 to Theorem 7.1.}  {\it If $3 \mid h(K)$ but $9 \nmid H_{-d}'(6)$, there is a nontrivial solution of the cubic Fermat equation in $K =\mathbb{Q}(\sqrt{-d})$.} \bigskip

It will also be convenient to express this corollary in terms of the polynomial $p_d(x)=m_d^\phi(x) m_d(x)$, where $\phi$ is the automorphism of $\Sigma/\mathbb{Q}$ that interchanges $\alpha$ and $\beta$.   Since $\wp_3 \mid \beta$ we have
$$m_d^\phi(x)=\prod_{\sigma \in G}{(x-\beta^\sigma)} \equiv x^{h(-d)}-Tr_{\Sigma/K}(\beta)x^{h(-d)-1} \hspace{.07 in} (\textrm{mod} \hspace{.05 in} \wp_3^2).$$ \smallskip
On the other hand, $\beta$ and $\sigma_1(\alpha) = \frac{3(\alpha+6)}{\alpha-3}$ are conjugates over $K$, so that
\begin{eqnarray*}
Tr_{\Sigma/K}(\beta) &=& Tr_{\Sigma/K}\left( \frac{3(\alpha+6)}{\alpha-3} \right) \equiv \sum_{\sigma \in G}{\left( \frac{3\alpha}{\alpha-3}\right)^\sigma} \equiv 3\sum_{\sigma \in G}{\left( \frac{1}{1-3/\alpha}\right)^\sigma}\\
&\equiv& 3\sum_{\sigma \in G}{\left( 1+\frac{3}{\alpha} \right)^\sigma} \equiv 3h(K)  \hspace{.07 in} (\textrm{mod} \ \wp_3^2),
\end{eqnarray*}
using the fact that $(\alpha,3)=\wp_3'$, and hence $\wp_3 \mid \frac{3}{\alpha}$.  Assuming $3 \mid h(K)$, this gives that $m_d^\phi(x) \equiv x^{h(-d)}$ (mod $\wp_3^2$), so that
$$p_d(x) \equiv x^{h(-d)} m_d(x)  \ (\textrm{mod} \ \wp_3^2), \quad 3 \mid h(-d).$$ \smallskip
We obtain the following criterion.  \bigskip

\noindent {\bf Proposition 7.2.} {\it When $3 \mid h(-d)$, $Tr_{\Sigma/K}(1/\alpha) \equiv 0$ (mod $\wp_3^2$) if and only if the coefficient of $x^{h(-d)+1}$ in $p_d(x)$ is divisible by $9$.} \bigskip

By the examples in \S6 we know that Aigner's conjecture holds for the 4 fields $K=\mathbb{Q}(\sqrt{-d})$ for which $-d \equiv 1$ (mod $3$) and $h(K)=3$.  We can also use Theorem 7.1 to verify the truth of Aigner's conjecture when $h(K)=6$ or $9$.  For all but one of the cases in the following proposition, it suffices to use the Corollary to Theorem 6.5 and the fact that $p_d(x) \equiv q_d(x)$ (mod $3$).  This congruence holds because the roots of $p_d(x)$ have the form $\alpha=3+\gamma^3$, as $\gamma$ runs through the roots of $q_d(x)$.\bigskip

\noindent {\bf Proposition 7.3.} {\it Aigner's conjecture is true whenever $K=\mathbb{Q}(\sqrt{-d})$ with $-d \equiv 1$ (mod 3) and $h(-d)=6$ or $9$. In the following 11 cases, namely,
\begin{eqnarray*}
h(-d)=6 &:& d=4 \cdot 26, 4 \cdot 29, 4 \cdot 38, 4\cdot 53, 515, 707;\\
h(-d)=9 &:& d=419, 491, 563, 1187, 2003;
\end{eqnarray*}
the point $Q_K$ is a nontrivial point on $E(K)$.}  \medskip

\noindent {\it Proof.} The quadratic fields $K$ for which  $-d \equiv 1$ (mod 3) and $h(K)=6$ are the fields with $d=104$, when the automorphism  $\displaystyle \tau=\left(\Sigma/K,\wp_3 \right)$ has order 3, and those with $d=4 \cdot 29, 4 \cdot 38, 4\cdot 53, 515,707$, when $\tau$ has order 6.  This may be verified by finding all the solutions of $3^6=x^2+dy^2$ or $(x^2+dy^2)/4$ and computing the corresponding class numbers $h(-d)$.  \medskip

For $d=104$, $p_{104}(x)$ is the $12$-th degree polynomial in the factorization of the iterated resultant $R_3(x)$ in Section 4.  Factoring $p_{104}(3+x^3)$ yields the polynomial
\begin{eqnarray*}
q_{104}(x)&=&x^{12}+4x^{11}+10x^{10}+16x^9+74x^8+136x^7+106x^6+408x^5\\
&&+666x^4+432x^3+810x^2+972x+729\\
& \equiv & x^6(x^3+2x+1)(x^3+x^2+2x+1) \ (\textrm{mod} \ 3).
\end{eqnarray*}
The coefficient of $x^7$ is $136 \equiv 1$ (mod $3$), so the Corollary of Theorem 6.5 implies that $Q_K$ is a nontrivial point in $E(K)$.  We list the minimal polynomial $q_d(x)$ of $\gamma$ for each of the remaining values of $d$, along with their factorizations (mod 3):
\begin{eqnarray*}
q_{116}(x)&=&x^{12}-2x^{11}+4x^{10}+34x^9-4x^8+14x^7+290x^6+42x^5-36x^4\\
&& +918x^3+324x^2-486x+729\\
&\equiv& x^6(x^6+x^5+x^4+x^3+2x^2+2x+2) \ (\textrm{mod} \ 3);
\end{eqnarray*}
\begin{eqnarray*}
q_{152}(x)&=&x^{12}+6x^{11}+9x^{10}+38x^9+114x^8+22x^7+601x^6+66x^5+1026x^4\\
&&+1026x^3+729x^2+1458x+729\\
&\equiv& x^6(x^6+2x^3+x+1) \ (\textrm{mod} \ 3);
\end{eqnarray*}
\begin{eqnarray*}
q_{212}(x)&=&x^{12}+12x^{11}+76x^{10}+254x^9+604x^8+1108x^7+1826x^6+3324x^5\\
&&+5436x^4+6858x^3+6156x^2+2916x+729\\
&\equiv& x^6(x^6+x^4+2x^3+x^2+x+2) \ (\textrm{mod} \ 3);
\end{eqnarray*}
\begin{eqnarray*}
q_{515}(x)&=&x^{12}+24x^{11}+179x^{10}+8x^9+1566x^8-1064x^7+7207x^6-3192x^5\\
&&+14094x^4+216x^3+14499x^2+5832x+729\\
&\equiv& x^6(x^6+2x^4+2x^3+x+1) \ (\textrm{mod} \ 3);
\end{eqnarray*}
\begin{eqnarray*}
q_{707}(x)&=&x^{12}-10x^{11}+462x^{10}+1196x^9+4146x^8+6974x^7+4582x^6+20922x^5\\
&&+37314x^4+32292x^3+37422x^2-2430x+729\\
&\equiv& x^6(x^6+2x^5+2x^3+2x+1) \ (\textrm{mod} \ 3).
\end{eqnarray*}

\noindent The assertion for these fields follows from the Corollary to Theorem 6.5 using the fact that each of the irreducible 6th degree polynomials over $\mathbb{F}_3$ listed in the above congruences has a nonzero coefficient of $x$.   \medskip

There are 5 fields $K$ with $-d \equiv 1$ (mod 3) and $h(K)=9$, corresponding to the values listed in the proposition.  In all of these cases the Frobenius automorphism $\tau$ has order 9, by the factorizations of $R_1(x)$ and $R_3(x)$ in \S4.  The polynomials $q_d(x)$ for the last four values of $d$ are:
\begin{eqnarray*}
q_{491}(x)&=& x^{18}-5x^{17}+161x^{16}+418x^{15}+1059x^{14}+3667x^{13}+10561x^{12}+17474x^{11}\\
&&+36518x^{10}+85772x^9+109554x^8+157266x^7+285147x^6+297027x^5\\
&&+257337x^4+304722x^3+352107x^2-32805x+19683\\
&\equiv& x^9(x^9+x^8+2x^7+x^6+x^4+x^3+2x^2+2x+2) \ (\textrm{mod} \ 3);
\end{eqnarray*}
\begin{eqnarray*}
q_{563}(x)&=& x^{18}-19x^{17}+188x^{16}+765x^{15}-1092x^{14}+1861x^{13}+17529x^{12}-10466x^{11}\\
&&-4240x^{10}+140654x^9-12720x^8-94194x^7+473283x^6+150741x^5\\
&&-265356x^4+557685x^3+411156x^2-124659x+19683\\
&\equiv& x^9(x^9+2x^8+2x^7+x^4+x^2+2x+2) \ (\textrm{mod} \ 3);
\end{eqnarray*}
\begin{eqnarray*}
q_{1187}(x)&=& x^{18}+51x^{17}+3388x^{16}+23875x^{15}+103588x^{14}+279691x^{13}+647729x^{12}\\
&&+1690194x^{11}+3278680x^{10}+5162354x^9+9836040x^8+15211746x^7\\
&&+17488683x^6+22654971x^5+25171884x^4+17404875x^3+7409556x^2\\
&&+334611x+19683\\
&\equiv& x^9(x^9+x^7+x^6+x^5+x^4+2x^3+x+2) \quad (\textrm{mod} \hspace{.05 in} 3);
\end{eqnarray*}
\begin{eqnarray*}
q_{2003}(x)&=& x^{18}-94x^{17}+32310x^{16}+350556x^{15}+2724866x^{14}+13517266x^{13}\\
&&+43159873x^{12}+ 106774252x^{11}+239739364x^{10}+464084648x^9\\
&&+719218092x^8+960968268x^7+1165316571x^6+1094898546x^5\\
&&+662142438x^4+255555324x^3+70661970x^2-616734x+19683\\
&\equiv& x^9(x^9+2x^8+2x^5+x^4+x^3+x^2+x+2) \ (\textrm{mod} \ 3).
\end{eqnarray*}

\noindent In each of these cases the assertion follows from the Corollary to Theorem 6.5.  For $d=419$ we have
\begin{eqnarray*}
q_{419}(x)&=&x^{18}+18x^{17}+66x^{16}-92x^{15}+1254x^{14}-1358x^{13}+4785x^{12}+4508x^{11}\\
&&-5844x^{10}+45656x^9-17532x^8+40572x^7+129195x^6-109998x^5\\
&&+304722x^4-67068x^3+144342x^2+118098x+19683,
\end{eqnarray*}
for which $q_{419}(x) \equiv x^9(x^9+x^6+x^4+2x^2+2)$ (mod $3$).
Here the hypothesis of Theorem 6.5 does not hold.  However, the polynomial $p_{419}(x)$ is
\begin{eqnarray*}
p_{419}(x)&=&x^{18}+1938x^{17}+1598844x^{16}-7296032x^{15}+210116832x^{14}-83424320x^{13}\\
&&+5572113408x^{12}-19699084288x^{11}+75000228864x^{10}-291034399744x^9\\
&&+1601000957952x^8-3440158470144x^7+8483079352320x^6\\
&&-59155454435328x^5+24881284988928x^4-229506313420800x^3\\
&&+706394569310208x^2+382092054626304x+2886917746065408,
\end{eqnarray*}
and the coefficient of $x^{10}$ is $75000228864 \equiv 6$ (mod $9$).  Alternatively, the minimal polynomial of $\alpha$ over $K=\mathbb{Q}(\sqrt{-419})$ is
$$m_{419}(x)=x^9+(969+39\sqrt{-419})x^8+(11292-888\sqrt{-419})x^7-(79156+2876\sqrt{-419})x^6$$
$$-(197304+1560\sqrt{-419})x^5+(-1282144+49408\sqrt{-419})x^4+(2913120+164448\sqrt{-419})x^3$$
$$+(4039168-4480\sqrt{-419})x^2+(26063616-609024\sqrt{-419})x-17042432-2489344\sqrt{-419}.$$ \smallskip
The coefficient of $x$ in $m_{419}(x)$ is
$$c=26063616-609024\sqrt{-419}=2^8 \cdot 3 \cdot (33937-793\sqrt{-419}),$$ \smallskip
which is not divisible by $\wp_3^2=(3, \frac{1+\sqrt{-419}}{2})^2$.  This shows that $\wp_3^2$ does not divide $Tr_{\Sigma/K}(1/\alpha)$ and therefore that $Q_K$ is nontrivial, by Theorem 7.1.  $\square$ \bigskip

A similar analysis has also been used to check the following assertion. \bigskip

\noindent {\bf Proposition 7.4.}  {\it Aigner's conjecture is true for the 24 fields $K=\mathbb{Q}(\sqrt{-d})$ with $-d \equiv 1$ (mod $3$) and class number $h(K)=12$.  These are the fields for which $d$ is one of the integers in the following list:}
\begin{eqnarray*}
ord(\tau)=6&:& 440, 680, 728, 1067, 1235, 1547, 1892, 1955, 2132, 2387, 2555, 2627, 2795, 2867\\
ord(\tau)=12&:&356, 731, 755, 932, 1208, 1355, 1763, 2468, 2723, 4907.
\end{eqnarray*}

For all of the discriminants in the above list except one, the truth of Aigner's conjecture is a consequence of Theorem 7.1 (or Theorem 6.5).  For only one of these integers, namely $d=2132=4 \cdot 533=4 \cdot 13 \cdot 41$, does the hypothesis of Theorem 7.1 fail.  We check Aigner's conjecture for this discriminant in the following computation. \bigskip

{\bf Example.} The discriminant $-d=-2132$ is the smallest discriminant I have found for which $K=\mathbb{Q}(\sqrt{-d})$ has class number divisible $h(K)$ by $3$ and $\wp_3^2$ divides $Tr_{\Sigma/K}(1/\alpha)$.  In this case the minimal poynomial of $\alpha$ over $K$ is
$$m_{2132}(x)=x^{12}+(9393228+147942\sqrt{-533})x^{11}+(17581542922-1139676838\sqrt{-533})x^{10}$$
$$+(4233420285756-214472221260\sqrt{-533})x^9+(28090931203668+4607212526412\sqrt{-533})x^8$$
$$+(-448665169157496-40619438690976\sqrt{-533})x^7$$
$$+(3238916409263024+167548974849520\sqrt{-533})x^6$$
$$+(-15707644756406928-209057497048512\sqrt{-533})x^5$$
$$+(42438813525646032+101094723966192\sqrt{-533})x^4$$
$$+(-111186869940745056-2651748110716320\sqrt{-533})x^3$$
$$+(403627934868140832+2285891562046368\sqrt{-533})x^2$$
$$+(-819117113722300800+26904445557929280\sqrt{-533})x$$
$$557549011339707200-46324546936236800\sqrt{-533}.\eqno{(7.4)}$$
The coefficient of $x$ is
$$c_1=-819117113722300800+26904445557929280\sqrt{-533}\eqno{(7.5)}$$
$$=2^6\cdot 3^2\cdot 5 \cdot 11\cdot (-25855969498810+849256488571\sqrt{-533}).$$
However, the point $Q_K$ is nontrivial, which can be seen by reducing $Q_K$ modulo the prime divisor $\wp_{569}$ of $p=569=6^2+533$ in $R_K$ for which $\sqrt{-533} \equiv -6$ (mod $\wp_{569}$). Then
$$m_{2132}(x) \equiv (x+565)(x+397)(x+74)(x+332)(x+344)(x+520)$$
$$\times (x+73)(x+336)(x+67)(x+94)(x+490)(x+286) \ (\textrm{mod} \ \wp_{569}).$$
The roots of this congruence correspond to solutions $(\alpha_i,\beta_i)$ of $Fer_3$ (mod $\wp_{569}$).  For example $(502, 180)$ is the solution of $Fer_3$ in $\mathbb{F}_{569}$ corresponding to the factor $(x+67)$.  Using (4.0), this corresponds in turn to the point $(18, 501)$ in $E(\mathbb{F}_{569})$.  An extended calculation shows that the sum of the $12$ points in $E(\mathbb{F}_{569})$ corresponding to the $12$ factors above is the point $(13,462) \equiv Q_K$ (mod $\wp_{569}$).  Hence, $Q_K \not \in \{O, (3,0), (3,9)\}$.  \bigskip\medskip

We can also see that $Q_K$ is nontrivial for $d=2132$ using the following criterion. \bigskip

\noindent {\bf Theorem 7.5.} {\it If $3 \mid h(K)$ and the solution $(\alpha,\beta)$ of $Fer_3$ in the Hilbert class field $\Sigma$ of $K$ satisfies
$$Tr_{\Sigma/K}\left(\frac{1}{\alpha}\right)+9Tr_{\Sigma/K}\left(\frac{1}{\alpha}\right)^2-9Tr_{\Sigma/K}\left(\frac{1}{\alpha \alpha^\tau}\right) \not \equiv 0 \ (mod \ \wp_3^3)$$
in $K$, then the point $Q_K \in E(K)$ is nontrivial. If $\wp_3^2$ divides $Tr_{\Sigma/K}\left(\frac{1}{\alpha}\right)$, this is equivalent to the condition}
$$\frac{1}{9}Tr_{\Sigma/K}\left(\frac{1}{\alpha}\right) \not \equiv Tr_{\Sigma/K}\left(\frac{1}{\alpha \alpha^\tau}\right) \ (mod \ \wp_3).\eqno{(7.6)}$$  \medskip

\noindent {\it Proof.} The congruence (7.3) implies by induction for integral elements $z_i \in \textsf{K}_3$ that
$$\sum_{i=1}^n{(z_i,w(z_i))}=(\zeta, w(\zeta)) \ \ \textrm{on} \ E',$$
where
$$\zeta \equiv \sum_{i=1}^n{z_i}+9\sum_{i \neq j}{z_i z_j^3} \ (\textrm{mod} \ 27).$$
Putting $n=h(K)$ and $z_i=\frac{1}{\alpha^{\sigma_i}}$ in this last congruence, as $\sigma_i$ runs over the elements of $G=\textrm{Gal}(\Sigma/K)$, gives that
\begin{eqnarray*}
z(Q_K) &\equiv& Tr\left(\frac{1}{\alpha}\right)+9\sum_{\sigma_i \neq \sigma_j}{\frac{1}{\alpha^{\sigma_i}} \left(\frac{1}{\alpha^{\sigma_j}}\right)^3} \ (\textrm{mod} \ \wp_3^3)\\
&\equiv& Tr\left(\frac{1}{\alpha}\right)+9\sum_{\sigma_i \neq \sigma_j}{\frac{1}{\alpha^{\sigma_i}} \frac{1}{\alpha^{\tau \sigma_j}}} \ (\textrm{mod} \ \wp_3^3)\\
&\equiv& Tr\left(\frac{1}{\alpha}\right)+9\sum_{j=1}^n{\frac{1}{\alpha^{\tau \sigma_j}}\sum_{\sigma_i \neq \sigma_j}{\frac{1}{\alpha^{\sigma_i}}}} \ (\textrm{mod} \ \wp_3^3)\\
&\equiv& Tr\left(\frac{1}{\alpha}\right)+9\sum_{j=1}^n{\frac{1}{\alpha^{\tau \sigma_j}}\left(Tr\left(\frac{1}{\alpha}\right)-\frac{1}{\alpha^{\sigma_j}}\right)} \ (\textrm{mod} \ \wp_3^3)\\
&\equiv& Tr\left(\frac{1}{\alpha}\right)+9Tr\left(\frac{1}{\alpha}\right)^2 -9\sum_{j=1}^n{\frac{1}{\alpha^{\tau \sigma_j}}\frac{1}{\alpha^{\sigma_j}}} \ (\textrm{mod} \ \wp_3^3)\\
&\equiv& Tr\left(\frac{1}{\alpha}\right)+9Tr\left(\frac{1}{\alpha}\right)^2 -9Tr\left(\frac{1}{\alpha^{\tau+1}}\right) \ (\textrm{mod} \ \wp_3^3).
\end{eqnarray*}
As in the proof of Theorem 7.1, if $Q_K=O$, then $z(Q_K)=0$, and the last congruence then contradicts our assumption.  $\square$ \bigskip

In the case $d=2132=4\cdot 533$ considered in the above example, we have $Tr_{\Sigma/K}\left(\frac{1}{\alpha}\right)=\frac{-c_1}{c_0}$, where $c_0$ and $c_1$ are the values in (7.4) and (7.5).  Moreover, $\frac{c_1}{9} \equiv c_0 \equiv 1$ (mod $\wp_3$), so
$$\frac{1}{9} Tr_{\Sigma/K}\left(\frac{1}{\alpha}\right) \equiv 2 \ (\textrm{mod} \ \wp_3).\eqno{(7.7)}$$
On the other hand, the minimal polynomial of $\displaystyle \rho=\alpha \alpha^\tau=\frac{3\alpha(\beta+6)}{\beta-3}$
over $K$ is
$$m_\rho(x)=x^{12}+(-1929300724-43134248\sqrt{-533})x^{11}$$
$$+(-18447676150143440+2680888726964336\sqrt{-533})x^{10}$$
$$+(902248097999755590848-97472186084833478528\sqrt{-533})x^9$$
$$+(29731434406711187538305920-420218764418100938219200\sqrt{-533})x^8$$
$$+(23938974253208159595881348288+1119867690172839485520702592\sqrt{-533})x^7$$
$$+(1406819705371218929952833452416-3851086926696195731132339712\sqrt{-533})x^6$$
$$+(-6234195043927405255639032653056-521801068020311548488682752512\sqrt{-533})x^5$$
$$+(36770079877907873033495144579072-748594576109130276089039049728\sqrt{-533})x^4$$
$$+(214837076703198452673508694118400+75601146639009852285065968025600\sqrt{-533})x^3$$
$$+(5303492224765691165830824722534400+402042353280849024747144430080000\sqrt{-533})x^2$$
$$+(-32803188406896794115454101127680000+89588612171678772118568616960000\sqrt{-533})x$$
$$+(-832937724789889206554236662638080000-51656410690117380065564763729920000\sqrt{-533}),$$
and we have
$$Tr_{\Sigma/K}\left(\frac{1}{\alpha \alpha^\tau}\right) =\frac{-a_1}{a_0},$$
where
$$a_1=-32803188406896794115454101127680000+89588612171678772118568616960000\sqrt{-533},$$
$$a_0=-832937724789889206554236662638080000-51656410690117380065564763729920000\sqrt{-533}.$$
Since $a_1 \equiv 0$ (mod $\wp_3$) and $a_0 \equiv 1$ (mod $\wp_3)$, this gives
$$Tr_{\Sigma/K}\left(\frac{1}{\alpha \alpha^\tau}\right) \equiv 0 \ (\textrm{mod} \ \wp_3),$$
and comparing with (7.7) we see that condition (7.6) is satisfied.  Thus, Theorem 7.5 implies that $Q_K$ is nontrivial in this case. \medskip

In conclusion, it it quite interesting that congruence conditions, such as those in Theorems 6.5, 7.1, and 7.5, can be used to guarantee the {\it existence} of solutions to the diophantine equation $Fer_3$ in imaginary quadratic fields.  In particular, Theorem 7.1 makes quite plausible that there should be a positive proportion of negative fundamental discriminants $d_K \equiv 1$ (mod $3$) with $3 \mid h(K)$ for which $Q_K$ is nontrivial.  Theorems 7.1 and 7.5 suggest that this proportion should be quite large, namely, that for at least $26/27$ of the discriminants in question, which is more than $96 \%$, $Q_K$ should be a nontrivial point.

\section{References}

\noindent [1] A. Aigner, Weitere Ergebnisse \"uber $x^3+y^3=z^3$ in quadratischen K\"orpern, Monatshefte f\"ur Mathematik 56 (1952), 240-252. \medskip

\noindent [2] A. Aigner, Ein zweiter Fall der Unm\"oglichkeit von $x^3+y^3=z^3$ in quadratischen K\"orpern mit durch 3 teilbarer Klassenzahl, Monatshefte f\"ur Mathematik 56 (1952), 335-338. \medskip

\noindent [3] A. Aigner, Die kubische Fermatgleichung in quadratischen K\"orpern, J. reine angew. Math. 195 (1955), 3-17. \medskip

\noindent [4] A. Aigner, Unm\"oglichkeitskernzahlen der kubischen Fermatgleichung mit Primfaktoren der Art $3n+1$, J. reine angew. Math. 195 (1955), 175-179. \medskip

\noindent [5] B.J. Birch, Heegner Points: The Beginnings, in {\it Heegner Points and Rankin $L$-series}, H. Darmon and S.-W. Zhang, eds., MSRI Publications 49, Cambridge University Press, 2004, pp. 1-10.  \medskip

\noindent [6] B. J. Birch and N. M. Stephens, Computation of Heegner points, ch. 1 of {\it Modular Forms}, R. A. Rankin, ed., Ellis Horwood Limited, Chichester, England, 1984, pp. 13-41. \medskip

\noindent [7] R. Bradshaw and W. Stein, Heegner points and the arithmetic of elliptic curves over ring class extensions, J. Number Theory 132 (2012), 1707-1719. \medskip

\noindent [8] D. Byeon, Class numbers of quadratic fields $\mathbb{Q}(\sqrt{D})$ and $\mathbb{Q}(\sqrt{tD})$, Proc. Amer. Math. Soc. 132 (2004), 3137-3140. \medskip

\noindent [9] H. Cohn, Iterated Ring Class Fields and the Icosahedron, Mathematische Annalen 255 (1981), 107-122.  \medskip

\noindent [10] David A.Cox, {\it Primes of the Form $x^2+ny^2$; Fermat, Class Field Theory, and Complex Multiplication}, John Wiley and Sons, 1989. \medskip

\noindent [11] J.E. Cremona, {\it Algorithms for Elliptic Modular Curves}, Cambridge University Press, 1992. \medskip

\noindent [12] M. Deuring, Die Typen der Multiplikatorenringe elliptischer Funktionenk\" orper, Abh. Math. Sem. Hamb. 14 (1941), 197-272.  \medskip

\noindent [13] M. Deuring, Die Anzahl der Typen von Maximalordnungen einer definiten Quaternionenalgebra mit primer Grundzahl, Jahresber. Deutsch. Math. Verein. 54 (1944), 24-41. \medskip

\noindent [14] M. Deuring, Teilbarkeitseigenschaften der singul\"aren Moduln der elliptischen Funktionen und die Diskriminante der Klassengleichung, Commentarii math. Helvetici 19 (1946), 74-82.  \medskip

\noindent [15] M. Deuring, Die Klassenk\"orper der komplexen Multiplikation, Enzyklop\"adie der math. Wissenschaften I2, 23 (1958), 1-60. \medskip

\noindent [16] W. Franz, Die Teilwerte der Weberschen Tau-Funktion, J. reine angew. Math. 173 (1935), 60-64. \medskip

\noindent [17] R. Fricke, {\it Lehrbuch der Algebra}, II, III, Vieweg, Braunschweig, 1928. \medskip

\noindent [18] R. Fueter, \"Uber kubische diophantische Gleichungen, Commentarii Math. Helvetici 2(1930), 69-89. \medskip

\noindent [19] B.H. Gross, Heegner points on $X_0(N)$, Chapter 5 in: {\it Modular Forms}, R.A. Rankin, ed., Ellis Horwood Limited and Halsted Press, Chichester, 1984, pp. 87-105. \medskip

\noindent [20] B.H. Gross and D.B. Zagier, Heegner points and derivatives of $L$-series, Invent. Math. 84 (1986), 225-320. \medskip

\noindent [21] H. Hasse, Neue Begr\"undung der komplexen Multiplikation. I. Einordnung in die allgemeine Klassenk\"orpertheorie, J. reine angew. Math. 157 (1927), 115-139; paper 33 in {\it Mathematische Abhandlungen}, Bd. 2, Walter de Gruyter, Berlin, 1975, pp. 3-27. \medskip

\noindent [22] H. Hasse, Zum Hauptidealsatz der komplexen Multiplikation, Monatshefte f\"ur Math. u. Physik 38 (1931), 315-322; paper 35 in {\it Mathematische Abhandlungen}, Bd. 2, Walter de Gruyter, Berlin, 1975, pp. 53-60. \medskip

\noindent [23] H. Hasse, {\it Vorlesungen \"uber Klassenk\"orpertheorie}, Physica-Verlag, W\"urzburg, 1933. \medskip

\noindent [24] H. Hasse, {\it Number Theory}, Springer, Berlin, 2002. \medskip

\noindent [25] M. Jones and J. Rouse, Solutions of the cubic Fermat equation in quadratic fields, Intl. J. Number Theory 9 (2013), 1579-1591. \medskip

\noindent [26] D. S. Kubert and S. Lang, {\it Modular Units}, Grundlehren der mathematischen Wissenschaften 244, Springer, New York, 1981. \medskip

\noindent [27] S. Lang, {\it Elliptic Functions}, Addison-Wesley, Reading, Mass., 1973. \medskip

\noindent [28] R. Lynch and P. Morton, The quartic Fermat equation in Hilbert class fields of imaginary quadratic fields, http://arxiv.org/abs/1410.3008, to appear in Intl. J. Number Theory. \medskip

\noindent [29] P. Morton, The cubic Fermat equation and complex multiplication on the Deuring normal form, Ramanujan Journal of Math. 25 (2011), 247-275.  \medskip

\noindent [30] P. Morton, Explicit congruences for class equations, Functiones et Approximatio Commentarii Mathematici 51 (2014), 77-110.  \medskip

\noindent [31] P. Morton, Supersingular parameters of the Deuring normal form, Ramanujan Journal of Math 33 (2014), 339-366.  \medskip

\noindent [32] P. Morton, Solutions of diophantine equations as periodic points of $p$-adic algebraic functions I, http://arxiv.org/abs/1410.4618, submitted.  \medskip

\noindent [33] J. Nakagawa and K. Horie, Elliptic curves with no rational points, Proc. Amer. Math. Soc. 104 (1988), 20-24.  \medskip

\noindent [34] M. Ram Murty, {\it Introduction to $p$-adic Analytic Number Theory}, AMS/IP Studies in Advanced Mathematics, vol. 27, American Math. Soc. and International Press, Providence, 2002. \medskip

\noindent [35] P. Ribenboim, {\it 13 Lectures on Fermat's Last Theorem}, Springer, New York, 1979. \medskip

\noindent [36] R. Schertz, {\it Complex Multiplication}, New Mathematical Monographs, vol. 15, Cambridge University Press, 2010. \medskip

\noindent [37] J.H. Silverman, {\it Advanced Topics in the Arithmetic of Elliptic Curves}, in: Graduate Texts in Mathematics, vol. 151, Springer, New York, 1994. \medskip

\noindent [38] J.H. Silverman, {\it The Arithmetic of Elliptic Curves}, 2nd edition, in: Graduate Texts in Mathematics, vol. 106, Springer, New York, 2009. \medskip

\noindent [39] H. Weber, {\it Lehrbuch der Algebra}, vol. III, Chelsea Publishing Co., New York, reprint of 1908 edition. \medskip

\noindent Department of Mathematical Sciences \smallskip

\noindent Indiana University - Purdue University at Indianapolis (IUPUI) \smallskip

\noindent 402 N. Blackford St., LD 270 \smallskip

\noindent Indianapolis, Indiana, 46202 \smallskip

\noindent {\it e-mail}: pmorton@math.iupui.edu

\end{document}